\documentclass[12pt,a4paper]{article}
\usepackage{amsfonts,amsmath}
\usepackage{latexsym}
\allowdisplaybreaks
\usepackage{times}
\usepackage{hyperref}
\title{$\SL(m,\C)$-equivariant and translation covariant continuous tensor valuations} 

\author{Judit Abardia-Ev\'equoz\footnote{Supported by DFG grant AB 584/1-2}, 
K\'aroly J. B\"or\"oczky\footnote{Supported by grants NKFIH ANN 121649, NKFIH K 109789 and
NKFIH K 116451}, M\'aty\'as Domokos\footnote{Supported by grants NKFIH K 119934}, D\'avid Kert\'esz\footnote{Supported by NKFIH grant ANN 121649}}

\newcommand{\proof}{\noindent{\it Proof: }}
\newcommand{\proofbox}{\mbox{ $\Box$}\\}
\newcommand{\R}{\mathbb{R}}
\newcommand{\C}{\mathbb{C}}

\newcommand{\N}{\mathbb{N}}
\newcommand{\T}{\mathbb{T}}

\newcommand{\V}{\mathbb{V}}

\newcommand{\nc}{\text{nc}}
\newcommand{\SL}{\mathrm{SL}}
\newcommand{\SU}{\mathrm{SU}}

\newcommand{\Val}{\mathrm{Val}}
\newcommand{\GL}{\mathrm{GL}}
\newcommand{\K}{\mathcal{K}}

\newtheorem{lemma}{Lemma}[section]
\newtheorem{theo}[lemma]{Theorem}

\newtheorem{defi}[lemma]{Definition}
\newtheorem{coro}[lemma]{Corollary}

\newtheorem{prop}[lemma]{Proposition}

\newtheorem{note}{Remark} 
\def\note#1{\ifvmode\leavevmode\fi\vadjust{\vbox to0pt{\vss
 \hbox to 0pt{\hskip\hsize\hskip1em
\vbox{\hsize3.5cm\small\raggedright\pretolerance10000
 \noindent #1\hfill}\hss}\vbox to8pt{\vfil}\vss}}}

\begin{document}

\maketitle

\begin{abstract}
The space of continuous, $\SL(m,\C)$-equivariant, $m\geq 2$, and translation covariant valuations taking values in the space of real symmetric tensors on $\C^m\cong\R^{2m}$ of rank $r\geq 0$ is completely described. The classification involves the moment tensor valuation for $r\geq 1$ and is analogous to the known classification of the corresponding tensor valuations that are $\SL(2m,\R)$-equivariant, although the method of proof cannot be adapted. 
\end{abstract}

\section{Introduction}

Let $n\geq 2$, let $\V$ be a vector space of real dimension $n$, and let $\mathcal{A}$ be an abelian semigroup. Denote by $\K(\V)$ the space of convex bodies in $\V$ (i.e., compact and convex sets in $\V$) equipped with the Hausdorff metric.  
An operator $Z:\K(\V)\to\mathcal{A}$ is called a \emph{valuation} if 
$$
Z(K\cup L)+Z(K\cap L)=Z(K)+Z(L)
$$
whenever $K,L\in\K(\V)$ satisfy that $K\cup L\in\K(\V)$. Here $`+`$ denotes the operation of the semigroup $\mathcal{A}$.

One of the principal aims in the theory of valuations is to obtain characterization results for known operators as the only valuations satisfying certain simple geometric and topological properties. Nowadays valuations taking values in different semigroups have been largely studied. 
The first classification theorem goes back to 1952, when Hadwiger proved that, for $\V=\R^n$, the linear combinations of intrinsic volumes are the only continuous real-valued valuations being invariant under rigid motions of $\R^n$ (see \cite{Had57}). 

Hadwiger's result can be generalized in different directions. For instance, we can change the group acting on $\K(\V)$ and classify the continuous real-valued valuations invariant under the action of some group (acting transitively on the unit sphere). This direction of study gave rise to the development of the theory of continuous and translation invariant real-valued valuations and has important consequences in integral geometry. We refer the reader to \cite{Ale01,Ale03,Ber09,Ber-Fu,bernig.fu.solanes,Ber-Sol,Fu06_2,wannerer14} and references therein for some results in this direction. 

Another important and more recent generalization of Hadwiger's theorem consists on changing the target space. For instance, valuations taking values in the space of convex bodies, in concave (or other spaces of) functions, etc. have been considered (see, e.g., \cite{Ale17+,boroczky.ludwig,colesanti.ludwig.mussnig,ludwig02,ludwig05,ludwig.survey,ludwig.survey2}). In these cases, the action of a group $G$ acting both on $\K(\V)$ and $\mathcal{A}$ is also considered and usually a characterization result for different groups $G$ and actions is studied. The related problem of tensor valuations on lattice polytopes is discussed in the pioneering paper of Ludwig and Silverstein \cite{LuS17}.

In this paper, we will focus on the study of tensor-valued valuations. We prove a Hadwiger-type theorem for the continuous, $\SL(m,\C)$-equivariant and translation covariant valuations taking values in the space of  real symmetric tensors of any given rank. Next we fix the notation to be used.

For $n\geq 2$, $r\in\N$ and an $n$-dimensional real vector space $\V$, we write $\T^r(\V)$ to denote the 
${n+r-1 \choose r}$-dimensional space of symmetric $r$-tensors of $\V$ over $\R$. 
In particular, $\T^0(\V)=\R$ and $\T^1(\V)=\V$. We write $S_r$ to denote the group of all permutations of $\{1,\ldots,r\}$. For $r\geq 2$, the \emph{symmetric tensor product} of $x_1,\ldots,x_r\in \V$ is defined by
$$
x_1\odot\ldots\odot x_r=\frac1{r!}\sum_{\sigma\in S_r}x_{\sigma(1)}\otimes\ldots\otimes x_{\sigma(r)}.
$$
We set $x^{r}=x\odot\ldots\odot x=x\otimes\ldots\otimes x$ for $x\in\V$.
In addition, the group $\GL(\V,\R)$ acts naturally on $\T^r(\V)$ as follows: For $\varphi\in{\rm GL}(\V,\R)$ the natural action on $\T^r(\V)$ is given by
$$
\varphi\cdot (x_1\odot\ldots \odot x_r)=\varphi x_1\odot\ldots \odot\varphi x_r
$$
for $x_1,\ldots,x_r\in\V$. We note that in this paper, tensor product is always over the reals even if the vector space has a complex structure, say possibly $\V=\C^m$ where $n=2m$, and $\T^r(\C^m)$ still means symmetric $r$-tensors over the reals.

Given an action of a closed subgroup $G\subset {\rm GL}(\V,\R)$ on $\V$, we say that a valuation $Z:\K(\V)\to\T^r(\V)$ is \emph{$G$-equivariant} if $Z(\varphi(K))=\varphi Z(K)$ holds for any $\varphi\in G$ and  $K\in \mathcal{K}(\V)$. If $r=0$, then $G$-equivariance is equivalent with $G$ invariance.

In the following, for $\V=\R^n$, we set $\mathcal{K}(\V)=\mathcal{K}^n$.
We say that a tensor valuation $Z:\,\mathcal{K}^n\to \T^r(\V)$ is \emph{translation covariant} if for every  
$K\in \mathcal{K}^n$, we have
\begin{equation}
\label{covariant-def}
Z(K+y)=\sum_{j=0}^r Z^{r-j}(K)\odot \frac{y^{j}}{j!}
\end{equation}
as a function of $y\in\R^n$ where each $Z^{r-j}$ is a tensor valuation of rank $r-j$ with $Z=Z^r$.
We observe that if $r=0$, then translation covariance is equivalent with translation invariance.
If $r>0$ and $Z$ is  $G$-equivariant for a closed subgroup $G\subset {\rm GL}(n,\R)$, then so is each $Z^{r-j}$. 

The reason  for the normalization in (\ref{covariant-def}) introduced by McMullen \cite{McM97} is that for $j=0,\ldots,r-1$, we have
\begin{equation}
\label{covariant-coefficients}
Z^{r-j}(K+y)=\sum_{m=0}^{r-j} Z^{r-j-m}(K)\odot \frac{y^{ m}}{m!},
\end{equation}
and hence $Z^{r-j}(K)$ is also a translation covariant valuation. 

For $r\geq 0$, a basic example of translation covariant tensor-valued valuation is the \emph{moment tensor valuation}
$$
M^r(K)=\frac1{r!}\int_Kx^{ r}\,dx,
$$
which is $\SL(n,\R)$-equivariant. 
For $K\in \mathcal{K}^n$, we write $V(K)$ to denote the volume of $K$,
and hence  for $y\in\R^n$, we have
\begin{equation}
\label{moment-tensor}
M^r(K+y)=\sum_{j=0}^r M^{r-j}(K)\odot \frac{y^{ j}}{j!}
\mbox{ \ where \ }M^0(K)=V(K).
\end{equation}

Haberl and Parapatits \cite{HaP} characterized the moment tensor valuation as continuous, $\SL(n,\R)$-equivariant, and translation covariant tensor valuation. More precisely, they 
characterized all measurable $\SL(n,\R)$-equivariant tensor valuations on polytopes containing the origin. As a special case of the main result of \cite{HaP}, we have the following.

\begin{theo}[Haberl, Parapatits]
\label{HaberlParapatits}
Let $n\geq 2$ and $r\geq 0$. An operator $Z:\,\mathcal{K}^n\to \T^r(\R^n)$ is an $\SL(n,\R)$-equivariant and translation covariant continuous valuation if and only if $Z=c\cdot M^{r}$ for a $c\in \R$, if  $r\geq 1$, and $Z=c_1+c_2V$ for $c_1,c_2\in \R$, if $r=0$.
\end{theo}

The main result in \cite{HaP} culminates a series of papers devoted to the study of tensor valuations that are affine-equivariant. The weakening of the continuity hypothesis to the measurability was an important aim after the results for upper semi-continuous valuations. We refer the reader to \cite{Ale99_1,bernig.hug,HaP1,HaP0,hug.schneider14,ludwig_02,ludwig_03,ludwig_11,ludwig_reitzner} for results in this direction and on tensor valuations.

\medskip
In this paper, we consider $\V=\C^m\cong\R^{2m}$ and $\SL(m,\C)$ acting on $\V$. We prove that the moment tensor valuation is again essentially the only $\SL(m,\C)$-equivariant and translation covariant tensor-valued valuation. More precisely, we prove the following result.

\begin{theo}
\label{sl(m,C)equi}
Let $m\geq 2$ and $r\geq 0$. An operator $Z:\,\mathcal{K}(\C^m)\to \T^r(\C^{m})$ is an $\SL(m,\C)$-equivariant and translation covariant continuous valuation if and only if $Z=c\, M^{r}$
for a $c\in \R$, if $r\geq 1$, and $Z=c_1+c_2V$ for $c_1,c_2\in \R$,  if $r=0$.
\end{theo}

We first notice that the case $r=0$ is not new since it can be obtained as a direct consequence of  characterization of the $\SU(m)$-invariant and translation invariant real-valued valuations  by Alesker \cite{Ale04a} if $m=2$ and 
by Bernig \cite{Ber09}) if $m\geq 3$.

\begin{theo}[Alesker, Bernig]
\label{Bernig}
Let $m\geq 2$. An operator $Z:\,\mathcal{K}(\C^m)\to \R$ is an $\SL(m,\C)$  and translation invariant continuous valuation if and only if $Z=c_1+c_2V$ for $c_1,c_2\in \R$.
\end{theo}

In Section~\ref{sec-real-valued}, we will provide a direct argument leading to Theorem~\ref{Bernig}, also with the aim to enlighten the general case $r\geq 1$ in Theorem~\ref{sl(m,C)equi}.

\smallskip
We notice that, if $n=2m$, $m\geq 2$, then $\SL(m,\C)$ is a closed subgroup of $\SL(n,\R)$ and hence Theorem~\ref{sl(m,C)equi} implies Theorem~\ref{HaberlParapatits} when the dimension $n$ of $\R^n$ is even.

With the tools used in this paper, to weaken the continuity hypothesis
in Theorem~\ref{sl(m,C)equi} to measurability or even upper-semicontinuity is, in the opinion of the authors, out of reach. Indeed, results from the theory of continuous and translation invariant valuations together with the fact that, in some contexts, continuity implies smoothness are heavily used, for instance, to differentiate some functions appearing on the proof of Theorem~\ref{sl(m,C)equi}. 
We also note that the method of the proof  of Theorem~\ref{HaberlParapatits} by Haberl and Parapatits \cite{HaP},
which lead to results under only measurability assumptions, does not seem to be
adaptable to Theorem~\ref{sl(m,C)equi}. One of the main ideas in \cite{HaP} is the use of double pyramids, which can be seen as a generalization of simplices. As the group $\SL(n,\R)$  acts transitively on the space of simplices, the study of the image of a fixed simplex suffices to know the image of every simplex. Since the group $\SL(m,\C)$ acts no longer transitively on the space of simplices in $\R^{2m}$ a similar argument does not seem to work for Theorem~\ref{sl(m,C)equi}. 

It was the paper Abardia and  Bernig \cite{Aba11} that first considered valuations intertwining $\SL(m,\C)$ by providing a generalization of the seminal characterization result for the projection body operator obtained by Ludwig \cite{ludwig05}.

\medskip
The paper is organized as follows: In Section~\ref{sec-conjectures}, we present the main steps for the proof of Theorem~\ref{sl(m,C)equi}, and reduce it to showing the non-existence of non-trivial even or odd, continuous, $\SL(m,\C)$-equivariant and translation invariant tensor-valued valuations, see Proposition~\ref{sl(m,C)invariant}.
Starting from Section~\ref{sec-trans-invariant}, the sole task of the paper is to prove Proposition~\ref{sl(m,C)invariant}.
Section~\ref{sec-trans-invariant} reviews the fundamental properties of translation invariant continuous valuations,
and Section~\ref{sec-subspaces} discusses real subspaces of $\C^m$.
 Theorem~\ref{Bernig} (the case $r=0$ of Theorem~\ref{sl(m,C)equi})
is proved in Section~\ref{sec-real-valued}.
For even valuations, the proof of Proposition~\ref{sl(m,C)invariant} is treated in Section~\ref{sec-Zeven}.  
In the case of odd valuations,
Proposition~\ref{sl(m,C)invariant} is verified in Section~\ref{sec-Zodd}.
In both cases, the section is divided into subsections according to the degree $j$ of homogeneity of the valuation. Putting together the result obtained for the different homogeneity degrees, the result in the odd and even cases follows, and Proposition~\ref{sl(m,C)invariant} is, in this way, proved (cf. page~\pageref{finalproof}).

\section{Proof of Theorem~\ref{sl(m,C)equi}}
\label{sec-conjectures}

In this section, we present the main ideas of the proof of Theorem~\ref{sl(m,C)equi}, and how to reduce it to Proposition~\ref{sl(m,C)invariant}. 

We start with the following fact for tensor-valued valuations, which was shown by McMullen \cite{McM97} if $s=r$ and by Alesker \cite{Ale99} for $s<r$.

\begin{theo}[McMullen, Alesker]
\label{maincoeff}
If $n\geq 2$, $r\geq 1$, $0\leq s\leq r$ and the valuations $Z:\,\mathcal{K}^{n}\to \T^r(\R^n)$ 
and $Z^{r-j}:\,\mathcal{K}^{n}\to \T^{r-j}(\R^n)$, $j=0,\ldots,s$ satisfy
$$
Z(K+y)=\sum_{j=0}^s Z^{r-j}(K)\odot \frac{y^{ j}}{j!}
$$
for $K\in \mathcal{K}^{n}$ and $y\in\R^n$, then $Z^{r-s}$ is translation invariant.
\end{theo}

The first new result that we need for the proof of Theorem~\ref{sl(m,C)equi} is the following.

\begin{prop}
\label{sl(m,C)mainterm}
If $m\geq 2$, $r\geq 1$ and $Z:\,\mathcal{K}(\C^m)\to \T^r(\C^m)$ is an $\SL(m,\C)$-equivariant 
translation covariant continuous valuation such that $Z^0\equiv c$
for a constant $c\in\R$; namely,
$$
Z(K+y)=c\cdot y^{ r}+\sum_{j=0}^{r-1} Z^{r-j}(K)\odot \frac{y^{ j}}{j!}
$$
for every $y\in\C^{m}$ and $K\in \mathcal{K}(\C^m)$, then $c=0$.
\end{prop}
\proof Let $v_1,\ldots,v_m$ be a complex basis of $\C^m$, and  let
$\V={\rm lin}_{\R}\{v_1,\ldots,v_m\}$. We observe that ${\rm SL}(\V,\R)\subset {\rm SL}(m,\C)$
is a closed subgroup, and 
 the action
of $\varphi\in {\rm SL}(\V,\R)$ on $\C^m$ is defined by
$\varphi(iv)=i\varphi(v)$ for $v\in\V$. For $\varrho=1,\ldots,r$, 
we consider the basis of $\T^\varrho(\C^m)$ induced by the real basis
$v_1,iv_1,\ldots,v_m,iv_m$ of $\C^m$.
The induced action of ${\rm SL}(\V,\R)$ on $\T^\varrho(\C^m)$ leaves
$\T^\varrho(\V)$ invariant, and $\T^\varrho(\V)$ has an ${\rm SL}(\V,\R)$-invariant direct complement subspace
spanned by the elements of the basis of $\T^\varrho(\C^m)$ containing at least one
of $iv_1,\ldots,iv_m$, which subspace in turn
is the kernel of a linear projection $\psi:\,\T^\varrho(\C^m)\to \T^\varrho(\V)$ 
commuting with the action of ${\rm SL}(\V,\R)$.

For $K\in \mathcal{K}(\V)$ and $j=0,\ldots,r-1$, we set $\widetilde{Z}^{r-j}(K)=\psi Z^{r-j}(K)$ and 
$\widetilde{Z}(K)=\psi Z(K)$. In particular, 
$\widetilde{Z}:\,\mathcal{K}(\V)\to \T^r(\V)$ is an $\SL(m,\R)$-equivariant, translation covariant and continuous valuation such that if $K\in \mathcal{K}(\V)$, then
$$
\widetilde{Z}(K+y)=c\cdot y^{ r}+\sum_{j=0}^{r-1} \widetilde{Z}^{r-j}(K)\odot \frac{y^{ j}}{j!}
$$
for $y\in\V$. 

On the other hand, Theorem~\ref{HaberlParapatits} and \eqref{moment-tensor} yield that 
$\widetilde{Z}^0=c_0V_m$  for a constant $c_0\in\R$ where $V_m$ is the $m$-dimensional volume on $\V$.
Therefore $c= c_0V_m(K)$ for all $K\in \mathcal{K}(\V)$, proving  that $c=0$.
\proofbox

The following statement is the main novel ingredient of the proof of 
Theorem~\ref{sl(m,C)equi}, and the rest of this paper will be devoted to its proof.

\begin{prop}
\label{sl(m,C)invariant}
If $m\geq 2$, $r\geq 1$ and $Z:\,\mathcal{K}^{2m}\to \T^r(\R^{2m})$ is an
odd or even $\SL(m,\C)$-equivariant and translation invariant continuous valuation, then $Z$ is constant zero.
\end{prop}

\noindent{\bf Proof of Theorem~\ref{sl(m,C)equi} based  on Propositions~\ref{sl(m,C)mainterm}
and \ref{sl(m,C)invariant}: }
Let $Z:\,\mathcal{K}^{2m}\to \T^r(\R^{2m})$, $r\geq 1$ be an $\SL(m,\C)$-equivariant and translation covariant continuous valuation, and hence
$$
Z(K+y)=\sum_{j=0}^r Z^{r-j}(K)\odot \frac{y^{ j}}{j!}, \quad y\in\R^{2m},
$$
where each $Z^{r-j}(K)$ is an $\SL(m,\C)$-equivariant and translation covariant tensor valuation of rank $r-j$,
$j=0,\ldots,r$. According to Theorem~\ref{Bernig}, $Z^0=c_1+c_2V$ for $c_1,c_2\in \R$. It follows from (\ref{moment-tensor}) that 
$\widetilde{Z}:=Z-c_2M^r$ 
is an $\SL(m,\C)$-equivariant and translation covariant tensor valuation of rank $r$, and
\begin{equation}
\label{moment-difference}
\widetilde{Z}(K+y)=c_1\cdot \frac{y^{ r}}{r!}+
\sum_{j=0}^{r-1}\widetilde{Z}^{r-j}(K)\odot \frac{y^{ j}}{j!}.
\end{equation}

We suppose that $\widetilde{Z}$ is not constant zero, and seek a contradiction.
First Proposition~\ref{sl(m,C)mainterm} yields $c_1=0$. Therefore there exists
a maximal $j\in\{0,\ldots,r-1\}$ such that 
$\widetilde{Z}^{r-j}$ is not constant zero. For $\varrho=r-j\geq 1$, we deduce from
Theorem~\ref{maincoeff} that the $\T^\varrho(\R^{2m})$-valued 
$\SL(m,\C)$-equivariant continuous valuation
$\widetilde{Z}^{\varrho}$ is  actually translation invariant.

Now we consider the $\SL(m,\C)$-equivariant and translation invariant
continuous $\T^\varrho(\R^{2m})$-valued valuations 
\begin{eqnarray*}
Z^+(K)&=&\frac12(\widetilde{Z}^{\varrho}(K)+\widetilde{Z}^{\varrho}(-K)),\\
Z^-(K)&=&\frac12(\widetilde{Z}^{\varrho}(K)-\widetilde{Z}^{\varrho}(-K)). 
\end{eqnarray*}
Since $-{\rm id}_{\R^{2m}}$ commutes with all elements of $\SL(m,\C)$, it follows that
$Z^+$ is even and $Z^-$ is odd.
Therefore 
 Proposition~\ref{sl(m,C)invariant} yields that 
 $Z^+$ and $Z^-$, and in turn $\widetilde{Z}^\varrho=Z^++Z^-$ is constant zero.
 This is absurd, thus $\widetilde{Z}$ is constant zero, completing
 the proof of Theorem~\ref{sl(m,C)equi}.
\proofbox

Therefore all we are left to prove is Proposition~\ref{sl(m,C)invariant}.

\section{Translation invariant continuous valuations}
\label{sec-trans-invariant}

Let $\V$ be a finite dimensional real vector space. In this section, we survey known properties for continuous and translation invariant valuations $Z:\,{\cal K}^n\to\V$ for $n\geq 2$. Our discussion is mostly based on Alesker \cite{Ale14}, and 
provide arguments 
using well-known ideas only when the statement we need has not been explicitly stated or proved before. 
We recall that ${\cal K}^n$ denotes the space of compact convex bodies in $\R^n$ and fix a real scalar product on $\R^n$. 
For general results in the theory of convex bodies and valuations, we refer, e.g., to the books~\cite{gardner.book06,gruber.book,Sch14}.

We write ${\rm Val}$ to denote the Fr\'echet space of continuous and translation invariant valuations $Z:\,\mathcal{K}^n\to \R$ (see Alesker \cite{Ale00} for a description of the Fr\'echet space). 
Hence the 
 Fr\'echet space of continuous and translation invariant valuations $Z:\,\mathcal{K}^n\to \V$ is ${\rm Val}\otimes \V$ (remember that tensor products are always over $\R$ in this paper).
We say that a valuation $Z:{\cal K}^n\to\V$ is \emph{homogeneous of degree $j$} or simply $j$-homogeneous if $Z(\lambda K)=\lambda^j Z(K)$ for every $\lambda\geq 0$ and $K\in{\cal K}^n$. We denote by $\Val_j\subset\Val$ the subset of $j$-homogeneous real-valued valuations.
Moreover, $Z:{\cal K}^n\to\V$ can be written uniquely in the form $Z=Z^++Z^-$ where $Z^+$ is even and $Z^-$ is odd; namely, 
$Z^+(-K)=Z^+(K)$ and $Z^-(-K)=-Z^-(K)$. $\Val^+$ (resp.~$\Val^-$) denote the subspace of even (resp.~odd) valuations in $\Val$.
A typical example of an even valuation with degree of homogeneity $j$ is the $j$th intrinsic volume $V_j$, which coincide with the $j$-dimensional Lebesgue measure on compact convex sets of dimension at most $j$.

We define an action of ${\rm GL}(n,\R)$ to ${\rm Val}\otimes\V$ as follows.
 \begin{defi}
Let ${\rm GL}(n,\R)$ act on the finite dimensional vector space $\V$ and denote this action by $\varphi\cdot v$, $\varphi\in{\rm GL}(n,\R)$, $v\in\V$. Then, the \emph{action of ${\rm GL}(n,\R)$ on ${\rm Val}\otimes \V$} is given by
\begin{equation}
\label{action-smoothness}
(\varphi Z)(K)=\varphi\cdot Z(\varphi^{-1}K),
\end{equation}
where $Z\in {\rm Val}\otimes \V$, $\varphi\in {\rm GL}(n,\R)$ and $K\in\mathcal{K}^n$.

In particular, if $\V=\T^r(\R^n)$, then
$$
(\varphi Z)(K)=\varphi\cdot Z(\varphi^{-1}K).
$$
\end{defi}

\begin{defi}
Let ${\rm GL}(n,\R)$ act on the finite dimensional vector space $\V$ and let $G\subset {\rm GL}(n,\R)$ be a closed subgroup. We say that a valuation $Z:\,\mathcal{K}^n\to \V$ is \emph{$G$-equivariant} if it is invariant under the action (\ref{action-smoothness}) over $G$.

We denote by $({\rm Val}\otimes \V)^G$ the Fr\'echet subspace of $G$-equivariant valuations.
\end{defi}

In this paper, $\V$ is always a finite dimensional real vector space. 
As stated in the introduction, our main focus is the case 
$$
\V=\T(\R^{2m})=\T(\C^{m})=\oplus_{r\geq 0}\T^r(\C^{m})
$$ 
and $G=\SL(m,\C)$ where $\T^r(\C^{m})$ is the real ${2m+r-1 \choose r}$-dimensional
space of $r$th symmetric tensor power of $\C^m$ over $\R$.

Another essential notion in the theory of valuations and in this paper is that of smoothness. 
\begin{defi}
We say that a valuation $Z:\,{\cal K}^n\to\V$ is \emph{smooth} if the action (\ref{action-smoothness}) defines a smooth map ${\rm GL}(n,\R)\to {\rm Val}\otimes\V$. Equivalently, $Z$ is smooth if $\varphi\mapsto Z\circ \varphi^{-1}$
is a smooth map ${\rm GL}(n,\R)\to {\rm Val}\otimes\V$. 
\end{defi}
In this paper, we use the terms smooth and $C^\infty$ interchangeably.
We write ${\rm Val}^\infty$  to denote the Fr\'echet subspace of smooth elements of ${\rm Val}$,
which is a dense subspace according to Alesker's Irreducibility Theorem (see \cite{Ale01}). 
\begin{theo}[Alesker's irreducibility theorem]
The natural representations of the group ${\rm GL}(n,\R)$ in $\Val_j^{+}$ and in $\Val_j^{-}$ are irreducible, i.e., there is no proper closed ${\rm GL}(n,\R)$ invariant subspace.
\end{theo}

The Fr\'echet subspace of smooth elements of ${\rm Val}\otimes\V$ is denoted by $({\rm Val}\otimes\V)^\infty$.  It follows from classical results in representation theory (see, e.g., \cite[p.~32]{Wal88}) that
\begin{equation}
\label{Valinfty}
({\rm Val}\otimes \V)^\infty={\rm Val}^\infty\otimes \V.
\end{equation}

If $Z:\,{\cal K}^n\to\R$ is invariant
under a closed subgroup $G\subset O(n)$ acting transitively on $S^{n-1}$, then  
\begin{description}
\item{(i)} $Z$ is smooth according to Alesker 
\cite{Ale04}; 
\item{(ii)} $Z$ is even according to Bernig \cite{Ber11}.
\end{description}

Normal cycles provide a natural way to represent smooth valuations.
If $K\in\mathcal{K}^n$, the \emph{normal cycle} of $K$ is defined as the set $\nc(K)\subset S\R^n:=\R^n\times S^{n-1}$ given by
$$
\nc(K)=\{(x,v)\in S\R^n\,:\,x\in K, \langle v,x-y\rangle\geq 0\,\forall y\in K\}.
$$
We say that an $(n-1)$-form $\omega\in\Omega^{n-1}(S\R^n)\otimes\V$ is translation invariant if it depends only on its components on $S^{n-1}$.

\begin{coro}
\label{dense}
$({\rm Val}\otimes \V)^\infty$ is a dense subspace of $\Val\otimes\V$. Moreover, the elements of 
$({\rm Val}\otimes \V)^{\infty}$ are given by integration over the normal cycle of a translation invariant form, i.e., if $Z\in({\rm Val}\otimes \V)^{\infty}$, then there exists a translation invariant $\omega\in\Omega^{n-1}(S\R^n)\otimes\V$
such that
$$
Z(K)=\int_{\nc(K)}\omega \mbox{ \ \ for $K\in\mathcal{K}^n$.}
$$
\end{coro}
\proof The first statement simply follows from (\ref{Valinfty}) and Alesker's Irreducibility Theorem \cite{Ale04}.

The second statement is proved by Alesker \cite{Ale06} if $\V=\R$, and hence it follows again by (\ref{Valinfty}).
\proofbox

For a closed subgroup $G\subset {\rm GL}(n,\R)$, we recall that $({\rm Val}\otimes \V)^G$ denotes the subspace of $G$-equivariant valuations in ${\rm Val}\otimes \V$.  Similarly to the real-valued case, as observed by Alesker and Bernig (private communication), any $Z\in({\rm Val}\otimes \V)^G$, with $G$ a closed subgroup $G\subset O(n)$ acting transitively on $S^{n-1}$, is smooth. For the convenience of the reader, we give a proof following the arguments of Corollary~3.3 in Fu \cite{Fu06} and Theorem~4.1 in Bernig \cite{Ber12}.

\begin{prop}[Alesker, Bernig]
\label{equivariantZsmooth}
Let a closed subgroup $G\subset O(n)$ act  transitively on $S^{n-1}$, and let ${\rm GL}(n,\R)$ act on a finite dimensional real vector space $\V$. Then, ${\rm dim}({\rm Val}\otimes \V)^G<\infty$ and $(\Val\otimes\V)^G\subset(\Val\otimes\V)^{\infty}$, that is, if a continuous translation invariant valuation $Z:\,{\cal K}^n\to\V$ is $G$-equivariant, then $Z$ is smooth.
\end{prop}
\proof For the argument, we fix a base point $e\in S^{n-1}$.

Let $Z\in({\rm Val}\otimes \V)^G$ and hence, $Z\in{\rm Val}\otimes \V$. From Corollary~\ref{dense}, 
it follows that there exists a convergent sequence $\{Z_{(m)}\}_{m\in\N}\subset ({\rm Val}\otimes \V)^{\infty}$ that converges to $Z$. For every $m\in\N$, define $\widetilde{Z}_{(m)}\in({\rm Val}\otimes \V)^{\infty,G}$ by 
$$
\widetilde{Z}_{(m)}(K)=\int_G\varphi Z_{(m)}(K)d\mu(\varphi)
$$
where $\mu$ is the probability Haar measure on $G$.
Since $Z$ is $G$-equivariant, it follows that the sequence $\{\widetilde{Z}_{(m)}\}$ also converges to $Z$.
In addition, each $\widetilde{Z}_{(m)}$ is a smooth and $G$-equivariant $\V$-valued valuation. We deduce from the second statement of Corollary~\ref{dense} that each $\widetilde{Z}_{(m)}$ is given by integrating an $(n-1)$-form $\omega_{(i)}\in\Omega^{n-1}(S\R^n)\otimes\V$ on the corresponding normal cycle. As $\widetilde{Z}_{(m)}$ is $G$ invariant, we can assume that $\omega_{(m)}$ is also $G$ invariant. Now we use that the group $G$ acts transitively on $S^{n-1}$, and hence, the form $\omega_{(m)}$ is determined by the knowledge of it in a single point (note that $\omega_{(m)}$ is translation invariant). In particular, it is enough to know 
$\omega_{(m)}|_{(o,e)}\in\Lambda^{n-1}(T_{(o,e)}S\R^n)\otimes\V$ where $T_{(o,e)}S\R^n$ stands for the tangent space at $(o,e)$  and $\Lambda^{n-1}(T_{(o,e)}S\R^n)\otimes\V$ is a finite dimensional vector space. 

As the sequence $\{\widetilde{Z}_{(m)}\}_{i\in\N}$, converges to $Z$, we have that the sequence 
$\{\omega_{(m)}|_{(o,e)}\}$
is bounded. Combining this fact with the finite dimensionality of $\Lambda^{n-1}(T_{(o,e)}S\R^n)\otimes\V$ yields that the sequence $\{\omega_{(m)}|_{(o,e)}\}$ has a limit point in $\Lambda^{n-1}(T_{(o,e)}S\R^n)\otimes\V$, therefore
the sequence $\{\omega_{(m)}\}$ has a translation invariant and $G$-equivariant limit point
$\omega\in\Omega^{n-1}(S\R^n)\otimes\V$. This $\omega$ gives rise to a smooth valuation $\psi$ in 
$({\rm Val}\otimes \V)^{G}$. However, by the convergence and definition of $\{\widetilde{Z}_{(m)}\}$, we have $\psi=Z$.

Since the finite dimensional vector space $\Lambda^{n-1}(T_{(o,e)}S\R^n)\otimes\V$ does not depend on the choice of $Z\in({\rm Val}\otimes \V)^G$, and the correspondence $Z\mapsto \omega|_{(o,e)}$ constructed above is clearly injective, we can identify $({\rm Val}\otimes \V)^G$ with a subspace of $\Lambda^{n-1}(T_{(o,e)}S\R^n)\otimes\V$, verifying that $({\rm Val}\otimes \V)^G$ is finite dimensional. \proofbox

In the following, we study the decomposition of $\Val$ and $\Val\otimes\V$ in terms of the degree of the homogeneity of the valuations and describe some of these spaces. 

McMullen \cite{McM77} proved the following useful polynomial behavior of certain valuations:

\begin{theo}[McMullen decomposition] 
\label{McMullendec}
Let $Z:{\cal K}^n\to\V$ be a continuous and translation invariant valuation, $K\in{\cal K}^n$, and $\lambda\geq 0$. Then, 
\begin{equation}
\label{Zpolynomial}
Z(\lambda K)=\sum_{j=0}^{n}\lambda^jZ_j(K)
\end{equation}
where $Z_j$ is a translation invariant continuous valuation homogeneous of degree $j$, $j=0,\ldots,n$
($Z_j(\lambda K)=\lambda^jZ_j(K)$ for $K\in{\cal K}^n$ and $\lambda\geq 0$).
In particular, 
$$
{\rm Val}=\oplus_{j=0}^n {\rm Val}_j,
$$
where ${\rm Val}_j$ denotes the Fr\'echet space of continuous and translation invariant valuations homogeneous of degree $j$, $j=0,\ldots,n$.
\end{theo}

For $G\subset {\rm GL}(n,\R)$ a closed subgroup, if $Z$ is $G$-equivariant, then the same holds for each $Z_j$.

\smallskip
Let us consider the coefficients occurring in (\ref{Zpolynomial}) for a continuous and translation invariant valuation 
$Z:\,{\cal K}^n\to\V$. We have that $Z_0$ is constant, and, as proved by 
Hadwiger \cite{Had57}, $Z_n$ is a constant multiple of the volume of $K$, that is, there exists $c\in\V$ such that
\begin{equation}
\label{Zn-hom}
Z_n(K)=c\cdot V(K)\mbox{ \ for $K\in {\cal K}^n$}.
\end{equation}
The valuation $Z_{n-1}$ can also be described. A direct extension of McMullen's representation result, proved in \cite{McM80}, gives us the following representation. 

\begin{theo}[McMullen] 
\label{Zn-1-homth}
Let $Z_{n-1}:{\cal K}^n\to\V$ be a continuous and translation invariant valuation homogeneous of degree $n-1$. Then, there exists a continuous $1$-homogeneous function $f:\,\R^n\to\V$ 
($f(\lambda x)=\lambda f(x)$ for $x\in\R^n$ and $\lambda\geq 0$) such that
\begin{equation}
\label{Zn-1-hom}
Z_{n-1}(K)=\int_{S^{n-1}}f\,dS_K\mbox{ \ for $K\in {\cal K}^n$,}
\end{equation}
where $S_K$ denotes the surface area measure of $K$ (see Schneider \cite{Sch14}).
Moreover, $f$ is unique up to a linear function. In other words, for continuous $1$-homogeneous functions $f,\tilde{f}:\,\R^n\to\V$, we have
\begin{equation}
\label{Zn-1-equality}
\int_{S^{n-1}}f\,dS_K=\int_{S^{n-1}}\tilde{f}\,dS_K\mbox{ \ for all $K\in {\cal K}^n$}
\end{equation}
if and only if $f-\tilde{f}$ is a linear function on $\R^n$.

In addition, $f$ is odd if $Z_{n-1}$ is odd. 
\end{theo}

We recall that if $h_C$ is the support function of a $C\in {\cal K}^n$, then
$$
\int_{S^{n-1}}h_C\,dS_K=nV(K,\ldots,K,C).
$$
Here $V(K,\ldots,K,C)$ denotes the mixed volume with $(n-1)$-times the convex body $K$ and once the convex body $C$ (see \cite[Section~5]{Sch14} for more information on mixed volumes). 
We note that if $\varphi$ is a volume preserving linear transformation, then
$$
\int_{S^{n-1}}h_{\varphi C}\,dS_{\varphi K}=nV(\varphi K,\ldots,\varphi K,\varphi C)=nV(K,\ldots,K,C)=\int_{S^{n-1}}h_C\,dS_K.
$$
Since any continuous $1$-homogeneous function $f:\,\R^n\to\V$
can be approximated by differences of support functions (see \cite[Lemma~1.7.8]{Sch14}), 
we deduce that if $\varphi\in {\rm GL}(n,\R)$ 
with $\det \varphi=\pm 1$ and $K\in {\cal K}^n$, then
\begin{equation}
\label{Zn-1-transfer}
\int_{S^{n-1}}f\,dS_{\varphi K}=\int_{S^{n-1}}f\circ \varphi^{-t}\,dS_{K}
\end{equation}
where $\varphi^{-t}$ stands for the inverse of the transpose of $\varphi$.

The following Proposition~\ref{Zn-1smoothfsmooth} is also observed by Alesker and Bernig (private communication). Below we provide an argument due to Alesker.

\begin{prop}[Alesker, Bernig]
\label{Zn-1smoothfsmooth}
Using the notation as above, 
$Z_{n-1}$ is smooth if and only if $f$ is smooth on $\R^n\backslash\{o\}$.
\end{prop}
\proof If $f$ is smooth, then readily the same holds for $Z$.

 We may assume that $\V=\R$. Let $C(S^{n-1})$ be the Banach space of continuous functions on 
$S^{n-1}$ with the $L_\infty$ norm, and let ${\rm Val}_{n-1}$ be the Fr\'echet space of $(n-1)$-homogeneous continuous translation invariant valuations  on ${\cal K}^n$. We write $C_0$ to denote the closed subspace of $C(S^{n-1})$ orthogonal to the $n$-dimensional subspace $L_0$ of $C(S^{n-1})$ linear maps
in terms of the $L_2$ scalar product of functions induced by the integral of their product; namely, $g\in C_0$ holds for $g\in C(S^{n-1})$ if and only if
\begin{equation}
\label{C0def}
\int_{S^{n-1}} g(u)\cdot u\,d{\cal H}^{n-1}(u)=o
\end{equation} 
where ${\cal H}^{n-1}$ denotes the $(n-1)$-dimensional Hausdorff measure.
Since ${\cal H}^{n-1}$ is invariant under ${\rm SO}(n)$, we observe that 
\begin{equation}
\label{phiginC0}
g\circ \varphi\in C_0 \mbox{ for any $g\in C_0$ and $\varphi\in {\rm SO}(n)$.}
\end{equation}

Let us consider the positive definite matrix
$$
M=\int_{S^{n-1}} u\otimes u\,d{\cal H}^{n-1}(u)=\int_{S^{n-1}} u\cdot u^t\,d{\cal H}^{n-1}(u).
$$
It follows that for any $\psi\in C(S^{n-1})$, there exists a unique $l\in L_0$ such that $\psi-l\in C_0$, namely,
$l(u)=\langle c_\psi,u\rangle$ for
$$
c_\psi=M^{-1}\int_{S^{n-1}} \psi(u)\cdot u\,d{\cal H}^{n-1}(u).
$$
Therefore (\ref{Zn-1-equality}) yields that the continuous linear map $\Omega:\, C_0\to {\rm Val}_{n-1}$ is bijective where
$$
\Omega(g)(K)=\int_{S^{n-1}}g \,dS_K\mbox{ \ for $g\in C_0$.}
$$
It follows from the open mapping theorem that the linear map $\Omega^{-1}$ is also continuous, therefore it is smooth.

Now for the smooth valuation $Z_{n-1}$, we consider $f=\Omega^{-1}(Z_{n-1})\in C_0$ that satisfies \eqref{Zn-1-hom}.
The map $F:\,{\rm SO}(n)\to {\rm Val}_{n-1}$ defined by
$$
F(\varphi)(K)=Z_{n-1}(\varphi^{-1}K)=
\int_{S^{n-1}}f \,dS_{\varphi^{-1}K}=\int_{S^{n-1}}f\circ \varphi^{-1}\,dS_{K}
\mbox{ \ for $\varphi\in {\rm SO}(n)$}
$$
is smooth (compare (\ref{Zn-1-transfer})), and hence $\Omega^{-1}\circ F$ satisfying
\begin{equation}
\label{fSOnsmooth}
\Omega^{-1}\circ F(\varphi)= f\circ \varphi^{-1} \mbox{ for $\varphi\in {\rm SO}(n)$ is smooth, as well.}
\end{equation}

Finally, since $f$ is $1$-homogeneous, it is enough to prove that the restriction of $f$ to $S^{n-1}$ is smooth. However, for orthogonal $u,v\in S^{n-1}$, the directional derivative of $f$ in the direction of $v$  at $u$ can be calculated using rotations around the $(n-2)$-dimensional linear subspace orthogonal to ${\rm lin}\{u,v\}$, showing that $f$ is $C^\infty$ as well.
\proofbox

We end this section with two useful results about the determination of $j$-homogenous valuations by knowing its value
 on some convex bodies. 

\begin{theo}[Schneider--Schuster \cite{ScS06}]
\label{theoZjj+1}
Let $j\in\{1,\dots,n-1\}$ and let $Z_j:\,{\cal K}^n\to\V$ be a continuous and translation invariant valuation homogeneous of degree $j$. Then,
\begin{equation}
\label{Zjj+1}
Z_j(K)=o\mbox{ \ for all $K\in {\cal K}^n$ if } Z_j(K)=o
\mbox{ \ for all $K\in {\cal K}^n$ with ${\rm dim}\,K=j+1$}.
\end{equation}
\end{theo}

For even valuations we have more information. Again let  $Z_j:\,{\cal K}^n\to\V$ be a continuous, $j$-homogeneous and translation invariant valuation for $j\in\{1,\ldots,n-1\}$.
 For any linear subspace $L$ of dimension $j$, Hadwiger's theorem (\ref{Zn-hom}) provides the existence of 
${\rm Kl}_{Z_{j}}(L)\in\V$ such that 
$$
Z_j(K)={\rm Kl}_{Z_j}(L)V_j(K)\text{ for }K\subset L.
$$ 
 The Grassmannian manifold ${\rm Gr}_j(\R^n)$ of linear subspaces of dimension $j$ of $\R^n$ is a smooth real algebraic 
subvariety of the real projective space over $\Lambda^j(\R^n)$. 
In this sense, the Klain function ${\rm Kl}_{Z_j}:{\rm Gr}_j(\R^n)\to\V$ is continuous, and it is smooth if $Z_j$ is smooth. We recall that for a compact topological space $X$, $C(X)$ denotes the normed space of continuous functions on $X$ with the maximum norm.

\begin{theo}[Klain's injectivity theorem \cite{Klain2000}]
\label{theoKlaininj}
The map ${\rm Kl}:\Val^+_j\to C({\rm Gr}_j(\R^n))$ is injective.

In particular, if $Z_j:{\rm K}^n\to\V$ is an even, continuous, translation invariant and $j$-homogenous valuation and there exists $c\in\V$ such that
\begin{equation}
\label{Klaininj}
\mbox{${\rm Kl}_{Z_j}(L)=c$ for any linear $j$-dimensional subspace $L$, then 
$Z_j=c\cdot V_j$.}
\end{equation} 
\end{theo}

In the following, we write $o$ to dente the origin in $\C^m=\R^{2m}$.

\section{Real vector subspaces of $\C^m$}\label{sec-subspaces}
In this section, we introduce the notation for linear subspaces in $\C^m$ and some properties of their bases.

We identify the complex vector space $\C^m$, of real dimension $2m$, with $\R^{2m}$ by using the bijection $\C^m\to\R^{2m}$ given by $$
(z_1,\dots,z_m)=(x_1+iy_1,\dots,x_m+iy_m)\mapsto(x_1,\dots,x_m,y_1,\dots,y_m).
$$

If $L\subset\C^m\cong\R^{2m}$ is a real vector subspace, then $\C L$ denotes the minimal complex linear subspace of 
$\C^m$ containing $L$. Hence ${\rm dim}_\C \C L$ is the maximal number of vectors in $L$ independent over $\C$. We say that a $j$-dimensional real subspace $L\subset\C^m\cong\R^{2m}$
is of {\it maximal complex rank} if
${\rm dim}_\C \C L=\min\{j,m\}$.

Next we describe a natural basis of a real subspace $L$ of $\C^m \cong\R^{2m}$. We observe that $\C^m$, $m\geq 2$, has a natural Hermitian inner product, whose real part is a scalar product on the underlying $\R^{2m}$.

\begin{lemma}
\label{real-complex-basis}
Let $L$ be a real vector subspace of $\R^{2m}=\C^m$ for $m\geq 2$ with ${\rm dim}_{\R}L=j\geq 1$, and let $d$ be the maximal number of vectors in $L$ independent over $\C$. Then, there exist $v_1,\ldots,v_d\in L$ independent over $\C$ such that
$v_1,\ldots,v_d$ is a real orthonormal basis of $L$, if $j=d$, and $v_1,\ldots,v_d,iv_1,\ldots,iv_{j-d}$ is a real orthonormal basis of $L$ if $j>d$.
\end{lemma}
\proof Let $U=L\cap iL$ be a complex subspace of $\R^{2m}=\C^m$ with $k={\rm dim}_\C U$, and let $W$ be the real orthogonal complement of $U$ inside $L$ with $t={\rm dim}_\R W$, and hence $j=2k+t$. If $k\geq 1$, then we choose a Hermitian basis 
$u_1,\ldots,u_k$ of $U$, and if $t\geq 1$, then we choose a real orthonormal basis $w_1,\ldots,w_t$ of $W$. We claim that if $t\geq 1$, then
\begin{equation}
\label{sumwl}
\alpha_1w_1+\ldots+\alpha_tw_t\in U\mbox{ \ for $\alpha_1,\ldots,\alpha_t\in\C$ yields
 $\alpha_1=\ldots=\alpha_t=0$.}
\end{equation}
We write $\beta_l={\rm Re}\,\alpha_l$ and $\gamma_l={\rm Im}\,\alpha_l$ for $l=1,\ldots,t$, and set $w=\gamma_1w_1+\ldots+\gamma_tw_t\in W$.
It follows from the condition in \eqref{sumwl} that $iw\in L$, and hence $iw\in L\cap iL=U$. However, $U$ is a complex subspace, thus $w\in U\cap W$. We conclude that $w=0$, and hence $\gamma_1=\ldots=\gamma_t=0$.
Therefore the condition in \eqref{sumwl} implies $\beta_1=\ldots=\beta_t=0$, proving \eqref{sumwl}.

 If $U=L$, and hence $d=k$ and $j=2k$, then we choose $v_l=u_l$ for $l=1,\ldots,d$. If $k=0$, or equivalently, $U=\{o\}$, then
$j=d=t$ by \eqref{sumwl}, and we choose $v_l=w_l$ for $l=1,\ldots,d$.
 If both $k\geq 1$ and $t\geq 1$, then $j=2k+t$ and $d=k+t$ by \eqref{sumwl},
we choose $v_l=u_l$ for $l=1,\ldots,k$ and $v_{k+l}=w_l$ for $l=1,\ldots,t$.
\proofbox

Now we show that for our purposes, we may assume that $d=\min\{j,m\}$ in Lemma~\ref{real-complex-basis}.

\begin{lemma}
\label{dense2}
If $m\geq 2$ and $j=1,\dots,2m-1$, then the subset of all $j$-dimensional real subspaces of maximal complex rank constitutes a dense subset of ${\rm Gr}_j(\R^n)$.
\end{lemma}
\proof
If $j=1$, then the statement readily holds, thus we assume $j>1$. 
Let $k=\min\{j,m\}$.
We call a $j$-dimensional real subspace $L$ of $\R^{2m}$ of \emph{lower complex rank} if ${\rm dim}_\C \C L<k$.

We recall that the Grassmannian manifold ${\rm Gr}_j(\R^{2m})$ of linear subspaces of dimension $j$ of 
$\R^{2m}$ is a connected smooth real algebraic 
subvariety of the real projective space over $\Lambda^j(\R^{2m})$, and in particular,
 locally it can be parametrized by the real wedge product of $j$ independent vectors over $\R$.  Now if an
 $L\in{\rm Gr}_j(\R^{2m})$ is represented
by $v_1\wedge \ldots \wedge v_j\in \Lambda^j(\R^{2m})$ for vectors $v_1,\ldots,v_j\in L$ independent over $\R$, then $L$ is of lower complex rank if and only if for any $1\leq i_1<\ldots<i_k\leq j$, the complex wedge product
$$
v_{i_1}\wedge_\C\ldots\wedge_\C v_{i_k}=0\in\Lambda^k(\C^m).
$$
Therefore real $j$-dimensional subspaces of lower complex rank form a real projective algebraic subvariety $X$ of  
${\rm Gr}_j(\R^{2m})$.
Since there exists some real $j$-dimensional subspace $L$ of maximal complex rank, and
${\rm Gr}_j(\R^{2m})$ is smooth and connected, the real dimension of $X$ is smaller than that of 
${\rm Gr}_j(\R^{2m})$. We conclude that $j$-dimensional subspaces of maximal complex rank
form a dense subset of ${\rm Gr}_j(\R^{2m})$.
\proofbox

According to Lemma~\ref{dense2}, the next corollary follows from
Klain's Injectivity Theorem~\ref{theoKlaininj} if the valuation $Z_j$ is even, and from
McMullen's Theorem~\ref{Zn-1-homth} and  Schneider's and Schuster's Theorem~\ref{theoZjj+1}
if the valuation $Z_j$ is odd.

\begin{coro}
\label{dense_sufficient}
For $m\geq 2$, $j=1,\dots,2m-1$
and finite dimensional real vector space $\V$,  let $Z_j:\K^{2m}\to\V$ be a continuous translation invariant valuation homogeneous of degree $j$. 
\begin{description}

\item[(i)] If $Z_j$ is even and ${\rm Kl}_{Z_j}(L)=0$ for every real subspace $L\in{\rm Gr}_j(\R^{2m})$
of maximal complex rank, then $Z_j$ is constant zero.

\item[(ii)] If $Z_j$ is odd and for every real subspace $L\in{\rm Gr}_{j+1}(\R^{2m})$ of maximal complex rank,
the continuous function $f$ on $L$ associated to the restriction of $Z_j$  to $L$ by \eqref{Zn-1-hom}
is linear, then $Z_j$ is constant zero.
\end{description}
\end{coro}

\section{Real valued $\SL(m,\C)$ and translation invariant continuous valuations}
\label{sec-real-valued}

In this section we give a direct proof of Theorem~\ref{Bernig}, basing on ideas in Abardia, Bernig \cite{Aba11}, Abardia \cite{Aba12,Aba15}. The main motivation to treat this particular case is that some of the main ideas to prove the general case (see Sections~\ref{sec-Zeven} and \ref{sec-Zodd}) are already contained in this section.

\smallskip
We first reduce the proof of Theorem~\ref{Bernig} by using McMullen's decomposition and Klain's injectivity theorem as follows:

Let $m\geq 2$ and let $Z:\,{\cal K}^{2m}\to\R$ be an $\SL(m,\C)$ and translation invariant continuous valuation. From the McMullen's decomposition, it follows that $Z=\sum_{j=0}^{2m} Z_j$ where each $Z_j$ is an $\SL(m,\C)$ and translation invariant continuous valuation homogeneous of degree $j$,
$j=0,\ldots,2m$. As we have described, $Z_0=c_1\chi$ for a constant $c_1\in\R$, and $Z_{2m}=c_2V$ 
for a constant $c_2\in\R$. 

Therefore we have to verify that $Z_j=0$ for $j=1,\ldots,2m-1$.
Since $Z_j$ is invariant under ${\rm SU}(m)\subset {\rm SL}(m,\C)$ that acts transitively on $S^{2m-1}$, we have that $Z_j$ is even from Bernig \cite{Ber09}. Thus Corollary~\ref{dense_sufficient}(i) applies, and Theorem~\ref{Bernig} follows if for each $j=1,\ldots,2m-1$,
\begin{equation}
\label{Zkzero}
{\rm Kl}_{Z_j}(L)=0\mbox{ \, for all $L\in{\rm Gr}_j(\R^{2m})$ of maximal complex rank. }
\end{equation}

\begin{lemma}
\label{d<m}
If $m\geq 2$, $j<m$, and $L$ is a $j$-dimensional real vector subspace of $\R^{2m}=\C^m$ of maximal complex rank, then ${\rm Kl}_{Z_j}(L)=0$.
\end{lemma}
\proof According to Lemma~\ref{real-complex-basis}, there exist $v_1,\ldots,v_j\in L$ independent over $\C$ such that $v_1,\ldots,v_j$ is a real basis of $L$.
We extend $v_1,\ldots,v_j$ to a complex basis $v_1,\ldots,v_m$ of $\C^m$. Since $j<m$, there exists
a $\varphi\in {\rm SL}(m,\C)$ such that $\varphi v_l=2v_l$ for $l=1,\ldots,j$. For the $j$-dimensional simplex $K$ with vertices $o,v_1,\ldots,v_j$, we have $V_j(K)>0$ and
$$
{\rm Kl}_{Z_j}(L) V_j(K)=Z_j(K)=Z_j(\varphi K)={\rm Kl}_{Z_j}(L) V_j(\varphi K)=2^j{\rm Kl}_{Z_j}(L) V_j(K),
$$
and hence ${\rm Kl}_{Z_j}(L)=0$. 
\proofbox

\noindent{\bf Proof of Theorem~\ref{Bernig}: }
According to \eqref{Zkzero} and Lemma~\ref{d<m}, we may assume that $Z_j$ is a $j$-homogeneous valuation with
$$
j\in\{m,\ldots,2m-1\}.
$$ 
Hence Lemma~\ref{real-complex-basis} yields that
there exists a complex basis  $v_1,\ldots,v_m$ of $\C^m=\R^{2m}$ such that writing $v_{l+m}=iv_l$
for $l=1,\ldots,m$, the vectors $v_1,\ldots,v_j\in L$  form a real basis of $L$.

Let $K\subset L$ be a $j$-dimensional  crosspolytope with vertices $\pm v_1,\ldots,\pm v_j$.
We claim that if $\psi\in {\rm GL}(m,\C)$ with $\det_{\C}\psi\in \R\backslash\{0\}$, then \begin{equation}
\label{dethom}
Z_j(\psi K)=\left|\mbox{$\det_{\C}\psi$}\right|^{\frac{j}m}Z_j(K).
\end{equation}
To prove (\ref{dethom}), first we assume that $\det_{\C}\psi>0$. In this case,
we set $D=\det_{\C}\psi$, and hence $\varphi=D^{\frac{-1}m}\psi\in{\rm SL}(m,\C)$ satisfies
$$
Z_j(\psi K)=Z_j(\varphi  D^{\frac{1}m}K)=Z_j(D^{\frac{1}m}K)=D^{\frac{j}m}Z_j(K),
$$
proving (\ref{dethom}) if $\det_{\C}\psi>0$.

If $\det_{\C}\psi<0$ in (\ref{dethom}), then we consider $\tilde{\psi}\in {\rm GL}(m,\C)$ defined by
$$
\tilde{\psi}(v_m)=-\psi(v_m)\mbox{ \ and \ }
\tilde{\psi}(v_l)=\psi(v_l)\mbox{ for $l=1,\ldots,m-1$. }
$$
It follows that $\det_{\C}\tilde{\psi}=|\det_{\C}\psi|$.  
Since the complex linear map $v_m\mapsto -v_m$ and 
$v_l\mapsto v_l$ for $l=1,\ldots,m-1$ leaves $K$ invariant, we have
$\tilde{\psi} K=\psi K$. Thus we deduce
$$
Z_j(\psi K)=Z_j(\tilde{\psi} K)=\left(\mbox{$\det_{\C}\tilde{\psi}$}\right)^{\frac{j}m}Z_j(K)=
\left|\mbox{$\det_{\C}\psi$}\right|^{\frac{j}m}Z_j(K),
$$
completing the proof of (\ref{dethom}).

\medskip

To finish the proof of Theorem~\ref{Bernig}, we distinguish two cases depending on whether $j>m$ or $j=m$.\\

\noindent{\bf Case  $m<j<2m$:}

For every $\lambda>0$, we define $\psi \in {\rm GL}(m,\C)$ by $\psi v_m=\lambda v_m$
and $\psi v_l=v_l$ for $l=1,\ldots,m-1$.
In particular, (\ref{dethom}) yields that
$$
Z_j(\psi K)= \lambda^{\frac{j}m}Z_j(K)=
 \lambda^{\frac{j}m}{\rm Kl}_{Z_j}(L) V_j(K).
$$
On the other hand, we observe that $\psi(iv_l)=iv_l$, $j=1,\ldots,j-m$, and hence $\psi$ maps $L$ into $L$. The real determinant of the restriction of $\psi$ to $L$ is 
$\lambda$. Thus
$$
Z_j(\psi K)={\rm Kl}_{Z_j}(L) V_j(\psi K)=\lambda{\rm Kl}_{Z_j}(L)V_j(K).
$$
We deduce that, for every $\lambda >0$, 
$$
(\lambda^{\frac{j}m}-\lambda){\rm Kl}_{Z_j}(L) V_j(K)=0.
$$
Hence, using $j>m$, we obtain ${\rm Kl}_{Z_j}(L)=0$.\\

\noindent{\bf Case $j=m$:}

 For $t\in (-\frac{\pi}2,\frac{\pi}2)$, let $K_t$ be the $m$-dimensional crosspolytope with vertices  
$\pm[(\sin t) v_1+(\cos t)iv_2],\pm v_2,\ldots,\pm v_m$. We consider the complex linear map $\psi_t$ defined by 
$\psi_t(v_1)=(\sin t) v_1+(\cos t)iv_2$
and $\psi_t(v_l)=v_l$ for $l=2,\ldots,m$. Thus 
$$
{\rm det}_{\C}\psi_t=\sin t.
$$
In addition, the $\varphi_t\in {\rm GL}(2m,\R)$ defined by 
$\varphi_t(v_1)=(\sin t) v_1+(\cos t)iv_2$, $\varphi_t(iv_2)=(-\cos t) v_1+(\sin t)iv_2$,
$\varphi_t(iv_1)=iv_1$, $\varphi_t(v_l)=v_l$ for $l\geq 2$
and $\varphi_t(iv_l)=iv_l$ for $l>2$ satisfies that
$K_t=\psi_t K=\varphi_t K$. 

We claim that
\begin{equation}
\label{Zgammatm}
Z_m(\varphi_t K)=Z_m(\psi_t K)=|\sin t|\cdot {\rm Kl}_{Z_m}(L)\cdot V_m(K).
\end{equation} 
Formula (\ref{Zgammatm}) follows from (\ref{dethom}) if $\sin t\neq 0$, and hence by the 
continuity of $Z_m$ if $\sin t=0$.

Now $Z_m$ is smooth because it is invariant under ${\rm SU}(m)$ (see Proposition~\ref{equivariantZsmooth}), and $\varphi_t\in {\rm GL}(2m,\R)$, $t\in (-\frac{\pi}2,\frac{\pi}2)$, is a $C^\infty$ family of $2m\times 2m$ matrices, thus  $Z_m(\varphi_t K)$ is a $C^\infty$ function of $t$. 
Since $Z_m(\varphi_t K)$ is differentiable at $t=0$, but the right-hand-side of \eqref{Zgammatm} is differentiable only if it vanishes, we conclude
${\rm Kl}_{Z_m}(L)=0$ by (\ref{Zgammatm}).

In turn, we deduce \eqref{Zkzero} for $j=m,\ldots,2m-1$.
Since Lemma~\ref{d<m} verifies \eqref{Zkzero}  for $j=1,\ldots,m-1$, 
the proof Theorem~\ref{Bernig} is now complete.
\proofbox

\section{$Z$ is even}
\label{sec-Zeven}

Let $r\geq 1$ and $m\geq 2$. For the whole section, fix an even, $\SL(m,\C)$-equivariant and translation invariant continuous valuation $Z:\,{\cal K}^{2m}\to\T^r(\R^{2m})$. By the McMullen decomposition (\ref{Zpolynomial}), we have $Z=\sum_{j=0}^{2m} Z_j$ where $Z_j$ is a $j$-homogeneous 
even $\SL(m,\C)$ and translation invariant continuous valuation for
$j=0,\ldots,2m$. 
We note that in this section, we do not use the inner product on $\R^{2m}$ at all. 

Proposition~\ref{sl(m,C)invariant} for even valuations will directly follow
after proving that the even valuation $Z_j$ is constant zero for each $0\leq j\leq 2m$, which we prove in the following. 

\medskip
Recall that ${\rm Gr}_j(\R^n)$ denotes the family of all real linear $j$-dimensional subspaces $L$ of $\R^{2m}$, 
$j=0,\ldots,2m$. For $j=0,\ldots,2m$ and $L\in{\rm Gr}_j(\R^n)$, we consider the Klain
constant ${\rm Kl}_{Z_j}(L)\in\T^r(\R^{2m})$ such that
$$
Z_j(K)={\rm Kl}_{Z_j}(L)V_j(K)\mbox{ \ for every $K\in\mathcal{K}^{2m}$ with $K\subset L$.}
$$
We recall that $V_j(K)$ is the $j$-dimensional volume of a compact convex $K\subset L$.

Since $Z_j$ is even and continuous, Klain's injectivity theorem (\ref{Klaininj}) applies, and Proposition~\ref{sl(m,C)invariant} for even valuations follows if
\begin{equation}
\label{Zjrevnzero}
{\rm Kl}_{Z_j}(L)=0\mbox{ \, for all $j=0,\ldots,2m$ and $L\in{\rm Gr}_j(\R^{2m})$}.
\end{equation}
More precisely, by Corollary~\ref{dense_sufficient}, we can reduce the problem to study only
  real $j$-planes of maximal complex rank in \eqref{Zjrevnzero}. Hence, to prove 
Theorem~\ref{sl(m,C)invariant} for even valuations, all we have to show
is that if $j\in\{0,\ldots,2m\}$ and $Z_j:\,{\cal K}^{2m}\to\T^r(\R^{2m})$ is a $j$-homogeneous 
even $\SL(m,\C)$-equivariant and translation invariant continuous valuation, then
\begin{equation}
\label{Zjrevnzero0}
{\rm Kl}_{Z_j}(L)=0\mbox{ \ for all  $L\in{\rm Gr}_j(\R^{2m})$ of maximal complex rank}.
\end{equation}
Hence, from now on, we always assume that
$$
\mbox{the  $L\in{\rm Gr}_j(\R^n)$ in \eqref{Zjrevnzero} is of maximal complex rank if $j=1,\dots,2m$.}
$$
According to Lemma~\ref{real-complex-basis}, there exists a complex basis $v_1,\ldots,v_m$ for $\R^{2m}=\C^m$
such that  setting $v_{m+l}=iv_l$ for $l=1,\ldots,m$, we have that $v_1,\ldots,v_m,v_{m+1},\ldots,v_{2m}$ form an $\R$-basis of $\R^{2m}$, and
$$
\mbox{ $v_1,\ldots,v_j$ form a real basis of  $L\in{\rm Gr}_j(\R^n)$}.
$$

We write $I$ to denote the family of all
 $\theta:\,\{1,\ldots,2m\}\to\N$ such that
\begin{equation}\label{defI}
\sum_{l=1,\ldots,2m}\theta(l)=r.
\end{equation}
It follows that ${\rm Kl}_{Z_j}(L)$ can be written in the form 
\begin{equation}
\label{Klexpansion}
{\rm Kl}_{Z_j}(L)=\sum_{\theta\in I}c_\theta\odot_{l=1}^{2m} v_l^{\theta(l)}
\end{equation}
where each $c_\theta:=c_{\theta,Z,j,L}\in\R$ depends on $\theta,Z,j,L$. 

For $j=0,\ldots,2m$ and $L\in{\rm Gr}_j(\R^n)$, let  $\psi\in {\rm SL}(m,\C)$ satisfy $\psi(L)= L$. Writing $\psi_L$ to denote the ($\R$-linear) restriction of $\psi$ to $L$, the core of our argument is the claim that
\begin{equation}
\label{dethomr}
\mbox{$|\det_{\R}\psi_L|$}\cdot{\rm Kl}_{Z_j}(L)=\psi\cdot{\rm Kl}_{Z_j}(L)
\end{equation}
where $|\det_{\R}\psi_L|=1$ if $j=0,m,2m$.
To prove \eqref{dethomr}, choose any $j$-dimensional compact convex set $K\subset L$, and hence
\begin{eqnarray*}
\left(\psi\cdot{\rm Kl}_{Z_j}(L)\right)V_j(K)&=&\psi\cdot Z_j(K)=Z_j(\psi K)={\rm Kl}_{Z_j}(L)V_j(\psi K)\\
&=&|\mathrm{det}_{\R}\psi_L|\cdot {\rm Kl}_{Z_j}(L)\,V_j(K).
\end{eqnarray*}
Now if $j=0$, then $V_j(\psi K)=V_j(K)=1$, thus $|\det_{\R}\psi_L|=1$.
If $j=2m$, then $\det_{\R}\psi_L=\det_{\R}\psi=|\det_{\C}\psi|^2=1$. Finally, if $j=m$, then
$\psi L=L$ yields that each entry of the matrix of $\psi\in {\rm SL}(m,\C)$ with respect to the complex basis $v_1,,\ldots,v_m$ of $\C^m$ is a real number, therefore
$\det_{\R}\psi_L=\det_{\C}\psi=1$.

We observe that if the map $\psi$ in (\ref{dethomr}) is the diagonal transformation with
$\psi(v_l)=\lambda_lv_l$  for $l=1,\ldots,m$
where each $\lambda_l>0$ and $\lambda_1\cdot\ldots \cdot\lambda_m=1$, then 
$\psi(v_{m+l})=\lambda_lv_{m+l}$ for $l=1,\ldots,m$ and
$\psi(L)=L$.  In this case, (\ref{dethomr}) is equivalent with the statement that for each $\theta\in I$, we have
\begin{equation}
\label{dethomrc}
\begin{array}{rcll}
c_\theta&=&
\left(\prod_{l=1}^m\lambda_l^{\theta(l)+\theta(m+l)}\right)\cdot c_\theta
&\mbox{ if $j=0,m,2m$};\\
\left(\prod_{l=1}^j\lambda_l\right)\cdot c_\theta&=&
\left(\prod_{l=1}^m\lambda_l^{\theta(l)+\theta(m+l)}\right)\cdot c_\theta
&\mbox{ if $j=1,\ldots,m-1$};\\
\left(\prod_{l=1}^{j-m}\lambda_l\right)\cdot c_\theta&=&
\left(\prod_{l=1}^m\lambda_l^{\theta(l)+\theta(m+l)}\right)\cdot c_\theta
&\mbox{ if $j=m+1,\ldots,2m-1$}.
\end{array}
\end{equation}
We also note that $\lambda_1\cdot\ldots \cdot\lambda_m=1$ yields
\begin{equation}
\label{dethomrc0}
\prod_{l=1}^m\lambda_l^{\theta(l)+\theta(m+l)}=
\prod_{l=1}^{m-1}\lambda_l^{\theta(l)+\theta(m+l)-\theta(m)-\theta(2m)}.
\end{equation}

Combining (\ref{dethomrc}) and (\ref{dethomrc0}), we deduce the following statements.

\begin{coro}
\label{Zjevenrj1m-1}
If $j=1,\ldots,m-1$ and $c_\theta\neq 0$ in (\ref{Klexpansion}), then
\begin{eqnarray*}
\theta(l)+\theta(m+l)&=&\theta(m)+\theta(2m)+1\mbox{ \ for $l=1,\ldots,j$,}\\
\theta(l)+\theta(m+l)&=&\theta(m)+\theta(2m)\mbox{ \ for $l=j+1,\ldots,m$}.
\end{eqnarray*}
In particular, $r=m(\theta(m)+\theta(2m))+j$.
\end{coro}

\begin{coro}
\label{ZjevenrjM+12m-1}
If $j=m+1,\ldots,2m-1$, $k=j-m$, and $c_\theta\neq 0$ in (\ref{Klexpansion}), then
\begin{eqnarray*}
\theta(l)+\theta(m+l)&=&\theta(m)+\theta(2m)+1\mbox{ \ for $l=1,\ldots,k$,}\\
\theta(l)+\theta(m+l)&=&\theta(m)+\theta(2m)\mbox{ \ for $l=k+1,\ldots,m$}.
\end{eqnarray*}
In particular, $r=m(\theta(m)+\theta(2m))+k$.
\end{coro}

In the following subsections, we prove that $Z_j\equiv 0$ for every $j=0,\dots,2m$ by distinguishing the different behaviors of $Z_j$ depending on $j$.

\subsection{Case $m+1\leq j\leq 2m-1$}
\label{sec-Zevenjmk}

\begin{lemma}
\label{m+1j2m-1}
$Z_{j}$ is constant zero for $j=m+1,\dots, 2m-1$.
\end{lemma}
\proof 
Let $k=j-m$, $k\in\{1,\ldots,m-1\}$. As in (\ref{Klexpansion}), we write
$$
{\rm Kl}_{Z_j}(L)=\sum_{\theta\in I}c_\theta\odot_{l=1}^{2m} v_l^{\theta(l)}.
$$
It follows from 
Corollary~\ref{ZjevenrjM+12m-1} that 
if $c_\theta\neq 0$ for $\theta\in I$, then 
$s:=\theta(m)+\theta(2m)$ satisfies that
$$
\theta(l)+\theta(m+l)=
\left\{ \begin{array}{rl}
s+1&\mbox{ \ for $l=1,\ldots,k$},\\
s &\mbox{ \ for $l=k+1,\ldots,m$}.
\end{array} \right.
$$
Let $\psi\in{\rm SL}(m,\C)$ be defined by $\psi(v_1)=-v_1$, $\psi(v_m)=-v_m$
and $\psi(v_l)=v_l$ if $1<l<m$, and hence
$$
\psi\cdot {\rm Kl}_{Z_j}(L)=(-1)^{(s+1)+s}{\rm Kl}_{Z_j}(L)=-{\rm Kl}_{Z_j}(L).
$$
This together with (\ref{dethomr}) implies ${\rm Kl}_{Z_j}(L)=0$, and in turn
we conlcude Lemma~\ref{m+1j2m-1} by Corollary~\ref{dense_sufficient} (i) for $j=m+1,\dots, 2m-1$.
\proofbox

\subsection{Case $1\leq j\leq m-1$}

\begin{lemma}
\label{Zj1m-1evenr}
$Z_j$ is constant zero for $j=1,\ldots,m-1$.
\end{lemma}
\proof
As in the proof of Lemma~\ref{m+1j2m-1}, we define $\psi\in\SL(m,\R)$ given by $\psi(v_1)=-v_1$, $\psi(v_{m})=-v_{m}$ and $\psi(v_l)=v_l$ for $2\leq l\leq m-1$. Applying \eqref{dethomr} to this $\psi$ and using the relations for $\theta$ given in Corollary~\ref{Zjevenrj1m-1}, we obtain that 
${\rm Kl}_{Z_j}(L)=-{\rm Kl}_{Z_j}(L)$, therefore ${\rm Kl}_{Z_j}(L)=0$. Using Corollary~\ref{dense_sufficient} (i), the statement of the lemma follows. 
\proofbox

\subsection{Case $j=m$}

In order to show that $Z_0$, $Z_m$  and $Z_{2m}$ are constant zero, we shall make use of the First Fundamental Theorem of classical invariant theory on $\SL(m,\R)$-invariants of several vectors. We describe it in the following. 

For $n\geq 2$, let $\V$ be an $n$-dimensional $\R$ vector space, and let $\T(\V)$ be the direct sum of all $\T^r(\V)$, $r\geq 0$. Hence $\T(\V)$ is an $\R$-algebra where the ``product'' is the symmetric tensor product. We observe that $\T(\V)$ can be naturally identified with the 
$\R$-algebra of polynomial functions on $\V^*$ where $\T^r(\V)$ corresponds to the homogeneous polynomials of degree $r$, and the identification respects the
${\rm GL}(\V,\R)$-action. 

For $m\geq 2$, we consider the standard representation of $\SL(m,\R)$ on the direct sum 
$\V=\R^m\oplus \R^m$. As the $\R$-algebras of symmetric tensors and polynomials can be identified,
we have the following consequence of the First Fundamental Theorem on vector invariants  of $\SL(m,\R)$ (see, e.g., Dolgachev \cite[Chapter~2]{Dol03}, Kraft, Procesi \cite[Section 8.4]{Kraft-Procesi} or Procesi \cite[Chapter 11.1.2]{procesi} for the general statement).

\begin{theo}[First Fundamental Theorem]
\label{fundamental}
Let $m\ge 2$, $r\geq 1$ and $\V=\R^m\oplus \R^m$, 
and let $\Theta\in\T^r(\V)$ be invariant under the natural action of  
 ${\rm SL}(m,\R)$.
\begin{description}
\item{(a)} If $m\geq 3$ or $r$ is odd, then $\Theta=0$;
\item{(b)} if $m=2$ and $r$ is even, then
$$
\Theta=c(v_1\odot w_2-v_2\odot w_1)^{r/2}
$$ 
where $c\in\R$ and $v_1,v_2$ form a basis of the first copy of $\R^2$, and $w_1,w_2$ 
is the corresponding basis of the second copy of $\R^2$.
\end{description}
\end{theo}

\begin{lemma}
\label{Zmevenr}
$Z_m$ is constant zero.
\end{lemma}
\proof According to Corollary~\ref{dense_sufficient}, it is sufficient to prove that
if $v_1,\ldots,v_m$ is a complex basis of $\C^m$, and $L={\rm lin}_\R\{v_1,\ldots,v_m\}$,
then ${\rm Kl}_{Z_m}(L)=0$. 
We observe that $\C^m=L\oplus_\R iL$ where  $iv_1,\ldots,iv_m$ is the corresponding real basis of $iL$.

It follows from (\ref{dethomr}) that
\begin{equation}
\label{dethomr2m}
{\rm Kl}_{Z_m}(L)=\psi\cdot{\rm Kl}_{Z_m}(L)
\end{equation}
for any $\psi\in {\rm SL}(m,\R)$. 
 We deduce from  the First Fundamental Theorem~\ref{fundamental} that ${\rm Kl}_{Z_m}(L)=0$ 
if $m\geq 3$ or $r$ is odd.

Therefore, we assume in the following that $m=2$ and $r$ is even. 
According to the First Fundamental Theorem~\ref{fundamental},
there exists $c\in\R$ such that writing  $w_1=iv_1$ and $w_2=iv_2$, we have
\begin{equation}
\label{Thetajm}
{\rm Kl}_{Z_2}(L)=c(v_1\odot w_2-v_2\odot w_1)^{r/2}.
\end{equation} 

We suppose that $c\neq 0$ in \eqref{Thetajm}, and seek a contradiction.
Let $K\subset L$ be the $2$-simplex with vertices $o,v_1,v_2$. For $t\in (-\frac{\pi}2,\frac{\pi}2)$, we define $K_t$ to be the $2$-dimensional simplex with vertices  $o,v_2,(\sin t) v_1+(\cos t)w_2$. 

\medskip
\noindent {\bf Claim 1:} \emph{If $c\neq 0$, then $r=2$, $Z(K_0)\neq 0$,
 and for any $t\in [0,\frac{\pi}2)$, we have}
\begin{equation}
\label{Zgammatmreven}
Z_2(K_t)=(\sin t) {\rm Kl}_{Z_2}(L) V_2(K)+(\cos t)Z_2(K_0).
\end{equation}

 For $t\in [0,\frac{\pi}2)$,
we consider the complex linear map $\varphi_t$ defined by 
$\varphi_t(v_1)=(\sin t) v_1+(\cos t)(iv_2)=(\sin t) v_1+(\cos t)w_2$
and $\varphi_t(v_2)=v_2$, thus 
$$
\mbox{$K_t=\varphi_tK$ and \ }{\rm det}_{\C}\varphi_t=\sin t.
$$
If $t\in (0,\frac{\pi}2)$, then $\psi_t=(\sin t)^{\frac{-1}2}\varphi_t\in {\rm SL}(2,\C)$  satisfies
\begin{eqnarray*}
\psi_t(w_1)&=&\psi_t(iv_1)=
(\sin t)^{\frac{-1}2}\Big((\sin t) w_1-(\cos t) v_2\Big)\\
\psi_t(w_2)&=&\psi_t(iv_2)=(\sin t)^{\frac{-1}2}\,w_2 .
\end{eqnarray*}
Since $Z_2$ is $2$-homogeneous, we deduce that if $t\in (0,\frac{\pi}2)$, then
\begin{eqnarray*}
Z_2(\varphi_t K)&=&Z_2((\sin t)^{\frac12}\psi_t K)=(\sin t)Z_2(\psi_t K)\\
&=&(\sin t)\psi_t\cdot Z_2(K)=
(\sin t)V_2(K)\psi_t\cdot{\rm Kl}_{Z_2}(L).
\end{eqnarray*}
For $t\in (0,\frac{\pi}2)$, we have
\begin{eqnarray*}
\psi_t\cdot{\rm Kl}_{Z_2}(L)&=&c(\sin t)^{\frac{-r}2}
\left((\sin t\, v_1+\cos t \,w_2)\odot w_2-v_2\odot(\sin t\, w_1-\cos t \,v_2)\right)^{\frac{r}2}\\
&=&c(\sin t)^{\frac{-r}2}
\left((\sin t (v_1\odot w_2-v_2\odot w_1)+\cos t (w_2\odot w_2+v_2\odot v_2)\right)^{\frac{r}2},
\end{eqnarray*}
implying the formula
\begin{eqnarray}
\label{rmightbe4}
Z_2(\varphi_t K)&=&cV_2(K)(\sin t)^{\frac{2-r}2}\\
\nonumber
&&\left((\sin t (v_1\odot w_2-v_2\odot w_1)+\cos t (w_2\odot w_2+v_2\odot v_2)\right)^{\frac{r}2}
\end{eqnarray}
Since $Z_2(K_t)=Z_2(\varphi_t K)$  is a continuous function of
$t\in [0,\frac{\pi}2)$, it follows that
\begin{equation}
\label{Z2K0limit}
\lim_{t\to 0^+}Z_2(\varphi_t K)=Z_2(\varphi_0 K)=Z_2(K_0).
\end{equation}
Combining  $r\geq 2$, \eqref{rmightbe4}, $cV_2(K)\neq 0$ and
$$
\lim_{t\to 0}\left(\sin t (v_1\odot w_2-v_2\odot w_1)+\cos t (w_2\odot w_2+v_2\odot v_2)\right)^{\frac{r}2}
=\left(w_2\odot w_2+v_2\odot v_2\right)^{\frac{r}2}\neq 0,
$$
we conclude that the limit in \eqref{Z2K0limit} exists only if $r=2$. Therefore $r=2$, and deduce from 
\eqref{rmightbe4} and \eqref{Z2K0limit} that
\begin{equation}
\label{Z2K0}
Z_2(K_0)=Z_2(\varphi_0 K)=cV_2(K)(w_2\odot w_2+v_2\odot v_2)\neq 0.
\end{equation}
We conclude  \eqref{Zgammatmreven}  if $t\in (0,\frac{\pi}2)$ from 
\eqref{Thetajm}, \eqref{rmightbe4} and \eqref{Z2K0}, and in turn if $t=0$ by continuity.\\

\noindent {\bf Claim 2:} \emph{If $c\neq 0$, then for any $t\in (-\frac{\pi}2,0)$, we have}
\begin{equation}
\label{Zgammatmrevenminus}
Z_2(\varphi_t K)=|\sin t| {\rm Kl}_{Z_2}(L) V_2(K)-(\cos t)Z_2(K_0).
\end{equation}

In this case, we have $\sin t<0$. The argument is similar as above only we modify the definition of $\varphi_t$ in order to have positive determinant and make use of the fact that we already know that $r=2$. For $t\in (-\frac{\pi}2,0)$, now  the complex linear map $\varphi_t$ is defined by  $\varphi_t(v_1)=v_2$,
$\varphi_t(v_2)=(\sin t) v_1+(\cos t)iv_2$. It follows that again $\varphi_t K=K_t=\varphi_t K$ and
$$
{\rm det}_{\C}\varphi_t=|\sin t|.
$$

Now $\psi_t=|\sin t|^{\frac{-1}2}\varphi_t\in {\rm SL}(2,\C)$ satisfies

\begin{eqnarray*}
\psi_t(w_1)&=&\psi_t(iv_1)=|\sin t|^{\frac{-1}2}\,w_2\\
\psi_t(w_2)&=&\psi_t(iv_2)= |\sin t|^{\frac{-1}2}
\Big( (\sin t) w_1-(\cos t)v_2\Big).
\end{eqnarray*}
Since $Z_2$ is $2$-homogeneous, we deduce
$$
Z_2(\varphi_t K)=|\sin t|Z_2(\psi_t K)=|\sin t|V_2(K)\psi_t\cdot{\rm Kl}_{Z_2}(L).
$$
As we already know that $r=2$ by Claim~1, in this case we have
\begin{eqnarray*}
\psi_t\cdot{\rm Kl}_{Z_2}(L)&=&c|\sin t|^{-1}
\left(v_2\odot(\sin t\, w_1-\cos t \,v_2)-(\sin t\, v_1+\cos t \,w_2)\odot w_2\right)\\
&=&c|\sin t|^{-1}
\left(-\sin t (v_1\odot w_2-v_2\odot w_1)-\cos t (w_2\odot w_2+v_2\odot v_2)\right),
\end{eqnarray*}
implying the formula
$$
Z_2(\varphi_t K)=cV_2(K)\left(|\sin t| (v_1\odot w_2-v_2\odot w_1)-\cos t (w_2\odot w_2+v_2\odot v_2)\right).
$$
In turn, we conclude \eqref{Zgammatmrevenminus} and Claim~2 if $t\in (0,\frac{\pi}2)$ by \eqref{Z2K0}.

It follows from  Claim~1, the continuity of $Z_2$ and Claim~2 that
$$
Z(K_0)=\lim_{t\to 0^+}Z_2(\varphi_t K)=\lim_{t\to 0^-}Z_2(\varphi_t K)=-Z(K_0),
$$
and hence $Z(K_0)=0$. This contradicts Claim~1, therefore proves ${\rm Kl}_{Z_2}(L)=0$
in \eqref{Thetajm}
for the case $m=2$ and $r$ is even, 
concluding the proof of Lemma~\ref{Zmevenr}.
\proofbox

\subsection{Case $j\in\{0,2m\}$}
\label{sec-Zevenj02m}

\begin{lemma}
\label{Z02mevenr}
$Z_0$ and $Z_{2m}$ are constant zero.
\end{lemma}
\proof Let $j\in\{0,2m\}$. According to (\ref{dethomr}), there exists a $\Theta\in\T^r(\R^{2m})$ such that
$Z_j(K)=\Theta V_j(K)$ for any $K\in\K(\R^{2m})$ and
\begin{equation}
\label{dethomr2m}
\Theta=\psi\cdot\Theta
\end{equation}
for any $\psi\in {\rm SL}(m,\C)$. 
In particular, we have that $\Theta\in\T^r(\R^m\oplus \R^m)$ is invariant under the natural action of ${\rm SL}(m,\R)$. We deduce from  the First Fundamental Theorem~\ref{fundamental} that $\Theta=0$ if $m\geq 3$ or $r$ is odd.

Therefore, we assume in the following that $m=2$ and $r$ is even. In this case, we choose a complex basis $v_1,v_2$ of $\C^2$, and define $w_l=iv_l$ for $l=1,2$. 
It follows from the First Fundamental Theorem~\ref{fundamental} that 
\begin{equation}
\label{Thetaj02m}
\Theta=c(v_1\odot w_2-v_2\odot w_1)^{r/2}
\end{equation} 
for a $c\in\R$. 
Since, by \eqref{dethomr2m}, $\Theta$ is not only invariant under $\SL(2,\R)$ but also under $\SL(2,\C)$, we consider 
$\psi\in {\rm SL}(2,\C)$ given by $\psi(v_1)=v_1$ and
$\psi(v_2)=iv_1+v_2=w_1+v_2$, and hence $\psi(w_1)=w_1$ and
$\psi(w_2)=-v_1+w_2$. 
A computation shows that
\begin{equation}
\label{psiThetaj2m}
\psi\cdot \Theta=c(v_1\odot w_2-v_2\odot w_1-v_1\odot w_1-v_2\odot w_2)^{r/2}.
\end{equation} 
If $c\neq 0$, then any term in $\Theta$ (see (\ref{Thetaj02m})) contains equal number of indices $1$ and $2$, while $\psi\cdot \Theta$ contains the term $(v_1\odot w_1)^{r/2}$ with non-zero coefficient (compare (\ref{psiThetaj2m})),
contradicting the invariance of $\Theta$ (see \eqref{dethomr2m}). Thus $c=0$, 
concluding the proof of Lemma~\ref{Z02mevenr}
\proofbox

\section{$Z$ is odd }
\label{sec-Zodd}

Let $m\geq 2$, $r\geq 1$, and let $Z:\,\mathcal{K}^{2m}\to \T^r(\R^{2m})$ be an odd $\SL(m,\C)$-equivariant and translation invariant continuous valuation, which we fix through the section. Similarly to Section 6, McMullen's decomposition theorem yields that
$Z=\sum_{j=0}^{2m}Z_j$ where each $Z_j$ is an odd $\SL(m,\C)$-equivariant and translation invariant continuous $j$-homogeneous valuation.
We prove in the following that $Z_j\equiv 0$ for every $0\leq j\leq 2m$. 

Let $j\in\{0,2m\}$. According to Hadwiger's theorem \eqref{Zn-hom}, there exists a constant $c_j\in\T^r(\R^{2m})$ such that
$Z_j(K)=c_j V_j(K)$ for any compact convex set $K$ in $\R^{2m}$. Since $Z_j$ is odd
and $V_j$ is even, we have
\begin{equation}
\label{Zoddj02m}
Z_j\equiv 0 \mbox{ \ if $j=0,2m$.}
\end{equation}

Therefore we may assume that $j\in\{1,\ldots,2m-1\}$. By Corollary~\ref{dense_sufficient} (ii),
it is sufficient to prove the following.

\begin{lemma}
\label{m-1to2m-1odd}
If $j=1,\ldots,2m-1$, $m\geq 2$, $Z_j:\,\mathcal{K}^{2m}\to \T^r(\R^{2m})$ is an odd $\SL(m,\C)$-equivariant and translation invariant continuous $j$-homogeneous valuation, $L\in{\rm Gr}_{j+1}(\R^{2m})$ is
of maximal complex rank, and the continuous 1-homogeneous function $f:\,L\to \T^r(\R^{2m})$ 
satisfies (cf. \eqref{Zn-1-hom}) 
\begin{equation}
\label{Zn-1-homL}
Z_j(K)=\int_{S^{2m-1}\cap L}f\,dS_{K,L}\mbox{ \ for $K\in {\cal K}(L)$,}
\end{equation}
where $S_{K,L}$ denotes the surface area measure  of $K$ with respect to $L$, then 
\begin{equation}
\label{alphabetaodd}
f(x+y)=f(x)+f(y)\mbox{ \ for any $x,y\in L$.}
\end{equation}
\end{lemma}

In order to prove that the function $f$ in Lemma~\ref{m-1to2m-1odd} is linear, we
 distinguish three different cases depending on whether $j\leq m-1$, $m\leq j\leq 2m-2$
or $j=2m-1$, but the idea described next is followed in all cases. 

It follows from Proposition~\ref{equivariantZsmooth} and $\mathrm{SU}(m)\subset\SL(m,\C)$ that $Z_j$ is smooth, and hence applying Proposition~\ref{Zn-1smoothfsmooth} to the restriction of
$Z_j$ to $\mathcal{K}(L)$ yields that
\begin{equation}
\label{fissmooth}
 \mbox{the function $f$ in \eqref{Zn-1-homL} is $C^\infty$ on $L\backslash\{o\}$}.
\end{equation}

As $L\in {\rm Gr}_{j+1}(\R^{2m})$ has maximal complex rank, 
Lemma~\ref{real-complex-basis} yields the existence of an orthonormal complex basis $v_1,\ldots,v_m$ for $\R^{2m}=\C^m$
such that  setting $v_{m+l}=iv_l$ for $l=1,\ldots,m$, we have that $v_1,\ldots,v_m,v_{m+1},\ldots,v_{2m}$ form an $\R$-basis of $\R^{2m}$, and
\begin{equation}
\label{v1vjL}
\mbox{ $v_1,\ldots,v_j$ form an orthonormal real basis of  $L\in{\rm Gr}_j(\R^n)$}.
\end{equation}
We note that the following ideas apply to any  complex basis $v_1,\ldots,v_m$ for $\R^{2m}=\C^m$
satisfying \eqref{v1vjL} where $v_{m+l}=iv_l$ for $l=1,\ldots,m$.

Similarly to the case of even valuations, we write $I$ to denote the family of all
 $\theta:\,\{1,\ldots,2m\}\to\N$ such that
\begin{equation}
\label{defIodd}
\sum_{l=1,\ldots,2m}\theta(l)=r.
\end{equation}
It follows that for any $x\in L$, $f(x)$ can be written in the form 
$$
f(x)=\sum_{\theta\in I}f_\theta(x)\odot_{l=1}^{2m} v_l^{\theta(l)}
$$
 where $f_\theta(x)\in\T^r(\R^{2m})$ for $x\in L$ and $\theta\in I$, and
\begin{equation}
\label{fthetaissmooth}
 \mbox{each $f_\theta$, $\theta\in I$, is $C^\infty$ on $L\backslash\{o\}$}
\end{equation}
according to \eqref{fissmooth}.

In order to prove~\eqref{alphabetaodd}, we use the ${\rm SL}(m,\C)$-equivariance of $Z_j$ as follows. Let $\varphi\in {\rm SL}(m,\C)$ satisfy that
$\varphi(L)=L$, and let $\Delta=\left|\det_\R\left(\varphi|_L\right)\right|$. It follows that
$\varphi|_L=\Delta^{\frac{1}{j+1}}\tilde{\varphi}$ 
where $\tilde{\varphi}\in{\rm GL}(L,\R)$ satisfies $\det_\R \tilde{\varphi}=\pm 1$.
Since 
$\varphi\cdot Z_j(K)=Z_j(\varphi K)$ for any $K\in \mathcal{K}(L)$, we deduce from 
\eqref{Zn-1-transfer} and \eqref{Zn-1-homL}  that
\begin{eqnarray*}
\int_{S^{2m-1}\cap L}\varphi\cdot f\,dS_{K,L}&=&
\int_{S^{2m-1}\cap L}f\,dS_{\varphi K,L}
=\int_{S^{2m-1}\cap L}f\,dS_{\Delta^{\frac{1}{j+1}}\tilde{\varphi}K,L}\\&=&
\int_{S^{2m-1}\cap L}\Delta^{\frac{j}{j+1}}f\,dS_{\tilde{\varphi}K,L}\\&=& 
\int_{S^{2m-1}\cap L}\Delta^{\frac{j}{j+1}}f\circ \tilde{\varphi}^{-t}\,dS_{K,L}.
\end{eqnarray*}
We conclude that
\begin{equation}
\label{phidefodd}
\Phi:=\varphi\cdot f-\left|\mathrm{det}_{\R} \left(\varphi|_{L}\right)\right|^{\frac{j}{j+1}}f\circ \tilde{\varphi}^{-t}
\end{equation}
is linear by (\ref{Zn-1-equality}). In particular, setting $k=j+1-m$ provided $j\geq m$, 
if $\alpha_1,\ldots,\alpha_m,\beta_1,\ldots,\beta_m\in\R$, then
\begin{equation}
\label{palphabetaodd}
\begin{array}{rlcl}
{\displaystyle\Phi\left(\sum_{q=1}^{j+1}\alpha_qv_q\right)-
\sum_{q=1}^{j+1} \Phi(\alpha_qv_q)}&=&0&\mbox{if $j\leq m-1$},\\
{\displaystyle\Phi\left(\sum_{q=1}^m\alpha_qv_q+\sum_{q=1}^k\beta_qiv_q\right)-
\sum_{q=1}^m \Phi(\alpha_qv_q)-\sum_{q=1}^k \Phi(\beta_qiv_q)}&=&0&
\mbox{if $j\geq m$}.
\end{array}
\end{equation}
The fact that \eqref{palphabetaodd} holds for some suitable family of possible $\varphi\in {\rm SL}(m,\C)$
will lead to \eqref{alphabetaodd}.

\subsection{Case $m\leq j\leq 2m-2$}

The whole section is devoted to prove the following statement. 

\begin{lemma}
\label{Zmleqjleq2m-2r}
If $j=m,\ldots,2m-2$, then
$Z_j$ is constant zero.
\end{lemma}

We prove Lemma~\ref{Zmleqjleq2m-2r} by a series of lemmas where we use 
the notation above set up  around Lemma~\ref{m-1to2m-1odd}. 
In particular, we fix a complex basis $v_1,\ldots,v_m$ of $\C^m$ such that
$v_1,\ldots,v_k,v_{k+1},\ldots,v_m,iv_1,\ldots,iv_k$ is a real orthonormal basis of its $\R$-linear span $L$, where $L$ is a $(j+1)$-dimensional real subspace with $j+1=k+m$.
For $\lambda>0$, $p\in\{1,\ldots,k\}$ and $l\in\{k+1,\ldots,m\}$, we frequently consider the map
$\varphi_{p,l}\in {\rm SL}(m,\C)$ defined by 
\begin{eqnarray*}
\varphi_{p,l}(v_p)&=&\lambda v_p, \\
\varphi_{p,l}(v_l)&=&\lambda^{-1} v_l, \\
\varphi_{p,l}(v_q)&=&v_q \mbox{ \ if $q\neq p,l$} 
\end{eqnarray*}
(we do not signal the dependence of $\varphi_{p,l}$ on $\lambda$).

In this case, $\varphi_{p,l}|_L$ is an $\R$-linear map of the $(k+m)$-dimensional subspace $L$ into $L$ whose determinant is $\lambda$. It follows that
$\varphi_{p,l}|_L=\lambda^{\frac{1}{m+k}}\tilde{\varphi}_{p,l}$ 
where $\tilde{\varphi}_{p,l}\in{\rm SL}(L,\R)$ satisfies 
\begin{eqnarray*}
\tilde{\varphi}_{p,l}(v_p)&=&\lambda^{\frac{j}{m+k}} v_p \mbox{ \ and \ } 
\tilde{\varphi}_{p,l}(v_{p+m})=\lambda^{\frac{j}{m+k}} v_{p+m},\\ 
\tilde{\varphi}_{p,l}(v_l)&=&\lambda^{-\frac{m+k+1}{m+k}} v_l,\\
\tilde{\varphi}_{p,l}(v_q)&=&\lambda^{-\frac{1}{m+k}} v_q 
\mbox{ \ if $q=1,\ldots,m+k$ and $q\neq p,l,p+m$},
\end{eqnarray*}
therefore $j=m+k-1$ implies
\begin{eqnarray*}
\tilde{\varphi}_{p,l}^{-t}(w_p)&=&\lambda^{-\frac{j}{m+k}} w_p,\\ 
\tilde{\varphi}_{p,l}^{-t}(w_l)&=&\lambda^{\frac{j+2}{m+k}} w_l,\\
\tilde{\varphi}_{p,l}^{-t}(w_q)&=&\lambda^{\frac{1}{m+k}} w_q 
\mbox{ \ if $q\neq p,l$.}
\end{eqnarray*}

\begin{lemma}
\label{insubspace}
If $w_q\in \C v_q\cap L$ for $q=1,\ldots,m$, then for any $\theta\in I$, we have
\begin{eqnarray}
\label{1tok}
f_\theta(w_1+\ldots+w_k)&=&f_\theta(w_1)+\ldots+f_\theta(w_k),\\
\label{k+1tom}
f_\theta(w_{k+1}+\ldots+w_m)&=&f_\theta(w_{k+1})+\ldots+f_\theta(w_m).
\end{eqnarray}
\end{lemma}
{\bf Remark } We observe that $w_q=\alpha_qv_q+\beta_qiv_q$ for $\alpha_q,\beta_q\in\R$ if $q=1,\ldots,k$,
and $w_q=\alpha_qv_q$ for $\alpha_q\in\R$ if $q=k+1,\ldots,m$.\\
\proof
To verify (\ref{1tok}) and (\ref{k+1tom}), it is sufficient to prove by induction for $p=1,\ldots,k$ and $l=k+1,\ldots,m$ that
if $w_q\in \C v_q\cap (L\backslash\{o\})$ for $q=1,\ldots,m$, then
\begin{eqnarray}
\label{1tokind}
f_\theta(w_1+\ldots+w_p)&=&f_\theta(w_1)+\ldots+f_\theta(w_p),\\
\label{k+1tomind}
f_\theta(w_{k+1}+\ldots+w_l)&=&f_\theta(w_{k+1})+\ldots+f_\theta(w_l)
\end{eqnarray}
by the continuity of $f_\theta$.

For \eqref{1tokind}, the case $p=1$ of the  induction argument trivially holds, therefore we assume that
$p>1$ and that \eqref{1tokind} holds for $p-1$.
We deduce from (\ref{palphabetaodd}) that for every $\theta\in I$, we have
\begin{eqnarray}
\label{fkmjdelta1top}
\lambda^{\theta(p)+\theta(p+m)-\theta(l)-\theta(l+m)}
\left(f_\theta(w_{1}+\ldots+w_p)-f_\theta(w_{1})-\ldots-f_\theta(w_p)\right)-&&\\
\nonumber
-\lambda^{\frac{j}{m+k}}
\left[f_\theta\left(\lambda^{\frac{-j}{m+k}}w_p+
\sum_{q=1}^{p-1}\lambda^{\frac{1}{m+k}}w_q
\right)-f_\theta\left(\lambda^{\frac{-j}{m+k}}w_p\right)-
\sum_{q=1}^{p-1}f_\theta\left(\lambda^{\frac{1}{m+k}}w_q
\right)
\right]
&=&0.
\end{eqnarray}
After multiplying \eqref{fkmjdelta1top} by $\lambda^{\theta(l)+\theta(l+m)-\theta(p)-\theta(p+m)}$, 
using the constant
$$
\delta_{p,l}=\theta(l)+\theta(l+m)-\theta(p)-\theta(p+m)+1,
$$
and using that $f(\alpha v)=\alpha f(v)$ for any $\alpha\in \R$ and $v\in L$, 
we deduce that (\ref{palphabetaodd}) is equivalent with
\begin{eqnarray}
\label{fdeltardelta1top}
f_\theta(w_{1}+\ldots+w_p)-f_\theta(w_{1})-\ldots-f_\theta(w_p)-&&\\
\nonumber
-\lambda^{\delta_{p,l}}
\left[f_\theta\left(\lambda^{-1} w_p+
\sum_{q=1}^{p-1}w_q
\right)-f_\theta(\lambda^{-1} w_p)\right]
+\lambda^{\delta_{p,l}}\sum_{q=1}^{p-1}f_\theta(w_q)
&=&0.
\end{eqnarray}
If $\delta_{p,l}<0$, then letting $\lambda$ tending to infinity in (\ref{fdeltardelta1top}), we deduce 
(\ref{1tokind}).

If $\delta_{p,l}>0$, then $f_\theta$ is differentiable at $w_p$ as $w_p\neq 0$ (compare \eqref{fissmooth}). In particular,
if $\lambda$ is small, then
\begin{eqnarray*}
f_\theta\left(\lambda^{-1} w_p+
\sum_{q=1}^{p-1}w_q
\right)-f_\theta(\lambda^{-1} w_l)
&=&\lambda^{-1}\left(f_\theta\left( w_p+
\lambda\sum_{q=1}^{p-1}w_q
\right)-f_\theta(w_p)\right)\\
&=&
\lambda^{-1}\,O\left( \lambda\sum_{q=1}^{p-1}w_q\right)=O(1),
\end{eqnarray*}
therefore letting $\lambda$ tend to zero in (\ref{fdeltardelta1top}) implies
(\ref{1tokind}).

Finally, if $\delta_{p,l}=0$, then (\ref{fdeltardelta1top}) reads
$$
f_\theta(w_{1}+\ldots+w_p)=(1-\lambda^{-1})f_\theta(w_{p})+f_\theta\left(\lambda^{-1} w_p+
\sum_{q=1}^{p-1}w_q
\right).
$$
Here letting $\lambda$ tending to infinity and applying the induction hypothesis, we complete the proof of
(\ref{1tokind}), and in turn (\ref{1tok}).

For \eqref{k+1tomind}, the case $l=k+1$  of induction argument trivially holds, therefore we assume that
$l>k+1$ and \eqref{k+1tomind} holds for $l-1$.
We deduce from (\ref{palphabetaodd}) that for every $\theta\in I$, we have
\begin{eqnarray}
\label{fkmjdeltak+1tom}
\lambda^{\theta(p)+\theta(p+m)-\theta(l)-\theta(l+m)}
\left(f_\theta(w_{k+1}+\ldots+w_m)-f_\theta(w_{k+1})-\ldots-f_\theta(w_m)\right)&&\\
\nonumber
-\lambda^{\frac{j}{m+k}}
\left[f_\theta\left(\lambda^{\frac{j+2}{m+k}}w_l+\lambda^{\frac{1}{m+k}}
\sum_{q=k+1}^{l-1}w_q
\right)-
f_\theta\left(\lambda^{\frac{j+2}{m+k}}w_l\right)-\lambda^{\frac{1}{m+k}}
\sum_{q=k+1}^{l-1}f_\theta(w_q)
\right]
&=&0.
\end{eqnarray}
After multiplying \eqref{fkmjdeltak+1tom} by $\lambda^{\theta(l)+\theta(l+m)-\theta(p)-\theta(p+m)}$, 
and using that $f(\alpha v)=\alpha f(v)$ for any $\alpha\in \R$ and $v\in L$, 
we deduce that (\ref{palphabetaodd}) is equivalent with
\begin{eqnarray}
\label{fdeltardeltak+1tom}
f_\theta(w_{k+1}+\ldots+w_l)-f_\theta(w_{k+1})-\ldots-f_\theta(w_l)-&&\\
\nonumber
-\lambda^{\delta_{p,l}}
\left[f_\theta\left(\lambda w_l+
\sum_{q=k+1}^{l-1}w_q
\right)-f_\theta(\lambda w_l)\right]
+\lambda^{\delta_{p,l}}\sum_{q=k+1}^{l-1}f_\theta(w_q)
&=&0.
\end{eqnarray}
If $\delta_{p,l}>0$, then letting $\lambda$ tending to zero in (\ref{fdeltardeltak+1tom}), we deduce 
(\ref{k+1tomind}).

If $\delta_{p,l}<0$, then $f_\theta$ is differentiable at $w_l$ as $w_l\neq 0$. In particular, if $\lambda$ is large,
then
\begin{eqnarray*}
f_\theta\left(\lambda w_l+
\sum_{q=k+1}^{l-1}w_q
\right)-f_\theta(\lambda w_l)&=&\lambda\left(f_\theta\left( w_l+
\frac{\sum_{q=k+1}^{l-1}w_q}\lambda
\right)-f_\theta( w_l)\right)\\
&=&
\lambda\,O\left( \frac{\sum_{q=k+1}^{l-1}w_q}\lambda\right)=O(1),
\end{eqnarray*}
therefore letting $\lambda$ tending to infinity in (\ref{fdeltardeltak+1tom}) implies
(\ref{k+1tomind}).

Finally, if $\delta_{p,l}=0$, then (\ref{fdeltardeltak+1tom}) reads
$$
f_\theta(w_{k+1}+\ldots+w_l)=(1-\lambda)f_\theta(w_{l})+f_\theta\left(\lambda w_l+
\sum_{q=k+1}^{l-1}w_q
\right).
$$
Here letting $\lambda$ tend to zero and applying the induction hypothesis complete the proof of
(\ref{k+1tomind}), and in turn (\ref{k+1tom}).
\proofbox

\begin{lemma}
\label{twosubspaces}
If $v\in{\rm lin}_\R\{v_1,iv_1,\ldots,v_k,iv_k\}$ and $w\in{\rm lin}_\R\{v_{k+1},\ldots,v_m\}$, then
$f_\theta(v+w)=f_\theta(v)+f_\theta(w)$.
\end{lemma}
\proof We may assume that $v,w\neq 0$. As $L\cap iL={\rm lin}_\R\{v_1,iv_1,\ldots,v_k,iv_k\}$ is a complex 
subspace of $\C^m$ of complex dimension $k$, we may may choose a complex basis $\tilde{v}_1,\ldots,\tilde{v}_k$
of $L\cap iL$ such that $v=\alpha_1\tilde{v}_1$ for $\alpha_1\in\R\backslash\{0\}$ and 
$\{\tilde{v}_1,i\tilde{v}_1,\ldots,\tilde{v}_k,i\tilde{v}_k\}$ is a real orthonormal basis. Similarly,
we may choose a real orthonormal basis $\tilde{v}_{k+1},\ldots,\tilde{v}_m$
of ${\rm lin}_\R\{v_{k+1},\ldots,v_m\}$ such that $w=\alpha_m\tilde{v}_m$ for $\alpha_m\in\R\backslash\{0\}$.
In particular, $\tilde{v}_1,\ldots,\tilde{v}_m$ is a complex basis of $\C^m$ such that 
$\tilde{v}_1,\ldots,\tilde{v}_m,i\tilde{v}_1,\ldots,i\tilde{v}_k$ form a real orthonormal basis of $L$.

For $\lambda>0$, 
we consider $\varphi\in {\rm SL}(m,\C)$ defined by $\varphi(\tilde{v}_1)=\lambda \tilde{v}_1$, 
$\varphi(\tilde{v}_m)=\lambda^{-1} \tilde{v}_m$ and $\varphi(\tilde{v}_q)=\tilde{v}_q$ if $q\neq 1,m$,
and hence the $\R$-linear map $\varphi_{1,m}|_L$ is of  determinant $\lambda$. Again, we do not signal the dependence of $\varphi$ on $\lambda$.
We have $\varphi|_L=\lambda^{\frac{1}{m+k}}\tilde{\varphi}$ 
where $\tilde{\varphi}\in{\rm SL}(L,\R)$ satisfies
\begin{eqnarray*}
\tilde{\varphi}^{-t}(\tilde{v}_1)&=&\lambda^{\frac{-j}{m+k}} \tilde{v}_1 \mbox{ \ and \ }
\tilde{\varphi}^{-t}(\tilde{v}_{m+1})=\lambda^{\frac{-j}{m+k}} \tilde{v}_{m+1} ,\\ 
\tilde{\varphi}^{-t}(\tilde{v}_m)&=&\lambda^{\frac{j+2}{m+k}} \tilde{v}_m\mbox{ \ and \ }
\tilde{\varphi}^{-t}(\tilde{v}_{2m})=\lambda^{-\frac{j}{m+k}} \tilde{v}_{2m},\\
\tilde{\varphi}^{-t}(\tilde{v}_q)&=&\lambda^{\frac{1}{m+k}} \tilde{v}_q 
\mbox{ \ if $q\neq 1,m,m+1,2m$},
\end{eqnarray*}
and hence $\tilde{\varphi}^{-t}(v)=\lambda^{\frac{-j}{m+k}} v$ and 
$\tilde{\varphi}^{-t}(w)=\lambda^{\frac{j+2}{m+k}}w$.

Based on the basis
$\tilde{v}_1,i\tilde{v}_1,\ldots,\tilde{v}_m,i\tilde{v}_m$ of $\R^{2m}$, $f$ can be written in the form 
$$
f(x)=\sum_{\theta\in I}\tilde{f}_\theta(x)\odot_{q=1}^{2m} \tilde{v}_q^{\theta(q)}
$$
for $x\in L$ and  $\tilde{f}_\theta(x)\in\R$. It is sufficient to prove that
$\tilde{f}_\theta(v+w)=\tilde{f}_\theta(v)+\tilde{f}_\theta(w)$ for $\theta\in I$.

 We deduce from the analogous of (\ref{palphabetaodd}) that for every $\theta\in I$, we have
\begin{eqnarray}
\label{fkmjdelta1m}
\lambda^{\theta(1)+\theta(1+m)-\theta(m)-\theta(2m)}
\left(\tilde{f}_\theta(v+w)-\tilde{f}_\theta(v)-\tilde{f}_\theta(w)\right)-&&\\
\nonumber
-\lambda^{\frac{j}{m+k}}
\left[\tilde{f}_\theta\left(\lambda^{\frac{-j}{m+k}}v+
\lambda^{\frac{j+2}{m+k}}w
\right)-\tilde{f}_\theta\left(\lambda^{\frac{-j}{m+k}}v\right)-
\tilde{f}_\theta\left(\lambda^{\frac{j+2}{m+k}}w
\right)
\right]
&=&0.
\end{eqnarray}

After multiplying \eqref{fkmjdelta1m} by $\lambda^{\theta(m)+\theta(2m)-\theta(1)-\theta(1+m)}$, 
using the constant
$$
\delta=\theta(m)+\theta(2m)-\theta(1)-\theta(1+m),
$$
and using that $\tilde{f}(\alpha v)=\alpha \tilde{f}(v)$ for any $\alpha\in \R$ and $v\in L$, 
we deduce that (\ref{palphabetaodd}) is equivalent with
\begin{eqnarray}
\label{fdeltardelta1m}
\tilde{f}_\theta(v+w)-\tilde{f}_\theta(v)-\tilde{f}_\theta(w)-&&\\
\nonumber
-\lambda^{\delta}
\left[\tilde{f}_\theta(v+\lambda^2w)-\lambda^2\tilde{f}_\theta(w)\right]
+\lambda^{\delta}\tilde{f}_\theta(v)
&=&0.
\end{eqnarray}
If $\delta>0$, then letting $\lambda$ tending to zero in (\ref{fdeltardelta1m}) yields
$\tilde{f}_\theta(v+w)=\tilde{f}_\theta(v)+\tilde{f}_\theta(w)$.

If $\delta<0$, then $\tilde{f}_\theta$ is differentiable at $w$ as $w\neq 0$, and hence
\begin{eqnarray*}
\tilde{f}_\theta(v+\lambda^2w)-\lambda^2\tilde{f}_\theta(w)&=& \lambda^2
\left(\tilde{f}_\theta\left(\frac{v}{\lambda^2}+w\right)-\tilde{f}_\theta(w)\right)\\
&=&
\lambda^2\,O\left(\frac{v}{\lambda^2} \right)=O(1)
\end{eqnarray*}
holds for large $\lambda$.
Therefore letting $\lambda$ tending to infinity in (\ref{fdeltardelta1m}) implies
$\tilde{f}_\theta(v+w)=\tilde{f}_\theta(v)+\tilde{f}_\theta(w)$.

Finally, if $\delta=0$, then (\ref{fdeltardelta1m}) reads
$$
\tilde{f}_\theta(v+w)=(1-\lambda^2)\tilde{f}_\theta(w)+\tilde{f}_\theta(v+\lambda^2w).
$$
Letting $\lambda$ tend to zero completes the proof of
$\tilde{f}_\theta(v+w)=\tilde{f}_\theta(v)+\tilde{f}_\theta(w)$.
\proofbox

Having Lemmas~\ref{insubspace} and \ref{twosubspaces} at hand, to prove that $f_\theta$ is linear,
 all we have to verify is that if $p=1,\ldots,k$, $\alpha_p,\beta_p\in\R\backslash\{0\}$ and $\theta\in I$, then
\begin{equation}
\label{alphabeta}
f_\theta(\alpha_p v_p+\beta_p iv_p)=f_\theta(\alpha_p v_p)+f_\theta(\beta_p iv_p).
\end{equation}

\begin{lemma}
\label{non-distribution}
If $\theta(p)+\theta(p+m)\neq \theta(m)+\theta(2m)$ for some $p\in\{1,\ldots,k\}$  and $\theta\in I$, then
\eqref{alphabeta} holds for all $\alpha_p,\beta_p\in\R$.
\end{lemma}
\proof  We may assume that $p=1$, and hence according to Lemma~\ref{twosubspaces} and the continuity of $f$, Lemma~\ref{non-distribution} is equivalent with the statement that if $\theta(1)+\theta(1+m)\neq \theta(m)+\theta(2m)$
for some $\theta\in I$, then
\begin{equation}
\label{f1lnon}
f_\theta(\alpha_1 v_1+\beta_1 iv_1+v_m)= f_\theta(\alpha_1 v_1)+f_\theta(\beta_1 iv_1)+f_\theta(v_m)
\end{equation}
for $\alpha_1,\beta_1\in\R\backslash \{0\}$.

It follows from (\ref{palphabetaodd}) that
\begin{eqnarray}
\label{fkmjdeltaplnon}
\lambda^{\theta(1)+\theta(1+m)-\theta(m)-\theta(2m)}
f_\theta(\alpha_1v_1+\beta_1iv_1+v_m)-&&\\
\nonumber
\lambda^{\theta(1)+\theta(1+m)-\theta(m)-\theta(2m)}
\left(f_\theta(\alpha_1v_1)+f_\theta(\beta_1iv_1)+f_\theta(v_m)\right)-&&\\
\nonumber
-\lambda^{\frac{j}{m+k}}
f_\theta\left(\lambda^{\frac{-j}{m+k}}(\alpha_1v_1+\beta_1iv_1)+
\lambda^{\frac{j+2}{m+k}}v_m\right)-&&\\
\nonumber
+\lambda^{\frac{j}{m+k}}
\left[f_\theta\left(\lambda^{\frac{-j}{m+k}}\alpha_1v_1\right)+
f_\theta\left(\lambda^{\frac{-j}{m+k}}\beta_1iv_1\right)
+f_\theta\left(\lambda^{\frac{j+2}{m+k}}v_m\right)
\right]
&=&0.
\end{eqnarray}
After multiplying \eqref{fkmjdeltaplnon} by 
$\lambda^{\theta(m)+\theta(2m)-\theta(1)-\theta(1+m)}$, 
using the constant
\begin{equation}
\label{deltapl0}
\delta=\theta(m)+\theta(2m)-\theta(1)-\theta(1+m)\neq 0,
\end{equation}
and using that $f(\alpha v)=\alpha f(v)$ for any $\alpha\in \R$ and $v\in L$, 
we deduce that (\ref{palphabetaodd}) is equivalent with
\begin{eqnarray}
\label{fdeltardeltaplnon}
f_\theta(\alpha_1v_1+\beta_1iv_1+v_m)-f_\theta(\alpha_1v_1)-f_\theta(\beta_1iv_1)-
f_\theta(v_m)-&&\\
\nonumber
-\lambda^{\delta}
\left[f_\theta\left(\alpha_1v_1+\beta_1iv_1+\lambda^{2}v_m\right)-f_\theta(\lambda^{2}v_m)\right]+&&\\
\nonumber
+\lambda^{\delta}(f_\theta(\alpha_1v_1)+f_\theta(\beta_1iv_1))
&=&0.
\end{eqnarray}

If $\delta>0$, then letting $\lambda$ tending to zero in (\ref{fdeltardeltaplnon}) implies
\begin{equation}
\label{f1lnon0}
f_\theta(\alpha_1v_1+\beta_1iv_1+v_m)-f_\theta(\alpha_1v_1)-f_\theta(\beta_1iv_1)-
f_\theta(v_m)=0,
\end{equation}
verifying \eqref{f1lnon}.

According to \eqref{deltapl0}, we may assume that $\delta<0$. In this case
$f_\theta$ is differentiable at $v_m$ by (\ref{fissmooth}), thus 
if $\lambda>0$ is small, then
\begin{eqnarray*}
f_\theta\left(\alpha_1v_1+\beta_1iv_1+\lambda^{2}v_m\right)-
f_\theta(\lambda^{2}v_m)&=&\lambda^{2}\left(
f_\theta\left( v_m+
\frac{\alpha_1v_1+\beta_1iv_1}{\lambda^{2}}
\right)-f_\theta( v_m)\right)\\
&=&
\lambda^{2}\,O\left(\frac{\alpha_1v_1+\beta_1iv_1}{\lambda^{2}}\right)=O(1),
\end{eqnarray*}
therefore letting $\lambda$ tending to infinity in (\ref{fdeltardeltaplnon}) leads to
\eqref{f1lnon0}, completing the proof of \eqref{f1lnon}. 
\proofbox

\noindent{\bf Proof of Lemma~\ref{Zmleqjleq2m-2r}: }
According to Lemmas~\ref{insubspace}, \ref{twosubspaces} and \ref{non-distribution},
to prove that the function $f$ in Lemma~\ref{m-1to2m-1odd} is linear,
 all we have to verify is that if $p\in\{1,\ldots,k\}$, $\alpha_p,\beta_p\in\R\backslash\{0\}$ and $\theta\in I$
satisfy $\theta(p)+\theta(p+m)=\theta(m)+\theta(2m)$, then
\begin{equation}
\label{alphabetafinal}
f_\theta(\alpha_p v_p+\beta_p iv_p)=f_\theta(\alpha_p v_p)+f_\theta(\beta_p iv_p).
\end{equation}

In order to prove the linearity for $f_\theta$, we define $\varphi\in\SL(m,\C)$ by $\varphi(v_p)=-v_p$, $\varphi(v_m)=-v_m$ and $\varphi(v_q)=v_q$ if $q\neq p,m$. 
In this case, $\varphi|_L$ is an $\R$-linear map of the $(m+k)$-dimensional subspace $L$ into $L$ whose determinant is $-1$. Moreover, $(\varphi|_{L})^{-t}=\varphi|_{L}$.
It follows from (\ref{palphabetaodd}) that
\begin{eqnarray*}
(-1)^{\theta(p)+\theta(p+m)+\theta(m)+\theta(2m)}
\left[f_\theta(\alpha_p v_p+\beta_p iv_p)-f_\theta(\alpha_p v_p)-f_\theta(\beta_p iv_p)\right]-&&\\
-\left[f_\theta(-\alpha_p v_p-\beta_p iv_p)-f_\theta(-\alpha_p v_p)-f_\theta(-\beta_p iv_p)\right]
&=&0
\end{eqnarray*}
Here $(-1)^{\theta(p)+\theta(p+m)+\theta(m)+\theta(2m)}=1$ as 
$\theta(p)+\theta(p+m)=\theta(m)+\theta(2m)$. Hence the fact that $f$ is odd yields
$$
2\left(f_\theta(\alpha_p v_p+\beta_p iv_p)-f_\theta(\alpha_p v_p)-f_\theta(\beta_p iv_p)\right)=0.
$$
We conclude \eqref{alphabetafinal}, and in turn Proposition~\ref{Zmleqjleq2m-2r}.
\proofbox

\subsection{Case $j=2m-1$}

The case $j=2m-1$ is handled similarly as the case $j=m,\ldots,2m-2$.

\begin{lemma}
\label{Z2m-1r}
$Z_{2m-1}$ is constant zero. 
\end{lemma}
\proof According to Lemma~\ref{m-1to2m-1odd}, it is sufficient to prove that the $f:\,\R^{2m}\to\T^r(\R^{2m})$ in Lemma~\ref{m-1to2m-1odd}
satisfies $f(x+y)=f(x)+f(y)$ for any
$x,y\in\C^m$. We claim that there exists a complex hermitian basis $v_1,\ldots,v_m$ of $\C^m$ such that
\begin{equation}
\label{uvv1v2}
x,y\in{\rm lin}_\R\{v_1,iv_1,v_2\}.
\end{equation}
Indeed, for \eqref{uvv1v2}, we may assume that $x$ and $y$ are complex independent. In particular, there exists
a complex hermitian basis $v_1,\ldots,v_m$ of $\C^m$ such that
 $x=\alpha\,v_1$ for $\alpha>0$, and the hermitian projection of $y$ onto the complex $(m-1)$-dimensional subspace of 
$\C^m$ complex orthogonal to $v_1$ is $\gamma v_2$ for $\gamma>0$, proving \eqref{uvv1v2}.

\smallskip
It follows from \eqref{uvv1v2} that it is sufficient to prove that if $v_1,\ldots,v_m$ is a complex hermitian basis of $\C^m$ and 
$\alpha_1,\beta_1,\alpha_2\in\R\backslash\{0\}$, then
\begin{equation}
\label{j2m-1flr}
f_\theta(\alpha_1v_1 +\beta_1iv_1+\alpha_2 v_2)=\alpha_1f_\theta(v_1) +\beta_1 f_\theta(iv_1)+\alpha_2 f_\theta(v_2),
\end{equation}
where we use the notation leading up to \eqref{palphabetaodd}.

For $\lambda>0$, we define $\varphi_\lambda\in {\rm SL}(m,\C)$ by $\varphi_\lambda(v_1)=\lambda v_1$, 
$\varphi_\lambda(v_2)=\lambda^{-1} v_2$
and $\varphi_\lambda(v_l)=v_l$ for $l>2$. It follows that
$$
\begin{array}{rlcrcl}
\varphi_\lambda^{-t}(v_1)=&\lambda^{-1}v_1 &\mbox{ \ and \ }
&\varphi_\lambda^{-t}(iv_1)&=&\lambda^{-1}iv_1,\\
\varphi_\lambda^{-t}(v_2)=&\lambda v_2 &\mbox{ \ and \ }&\varphi_\lambda^{-t}(iv_2)&=&\lambda v_2, \\
\varphi_\lambda^{-t}(v_l)=&v_l &\mbox{ \ and \ }&\varphi_\lambda^{-t}(iv_l)&=&iv_l \mbox{ \ for $l>2$.}
\end{array}
$$
We observe that for $\varphi=\varphi_\lambda$, we have
$\varphi|_L=\tilde{\varphi}=\varphi_\lambda$  in (\ref{palphabetaodd}).
Writing 
$$
\delta=\theta(1)+\theta(m+1)-\theta(2)-\theta(m+2)
$$
for $\theta\in I$,  (\ref{palphabetaodd}) yields that 
\begin{eqnarray*}
\lambda^{\delta} \left[f_\theta(\alpha_1v_1 +\beta_1iv_1+\alpha_2 v_2)-
 \alpha_1f_\theta(v_1)- \beta_1f_\theta(iv_1)-
 \alpha_2f_\theta(v_2)\right]-&&\\
- f_\theta(\lambda^{-1}\alpha_1v_1 +\lambda^{-1}\beta_1iv_1+\lambda\alpha_2 v_2)+&&\\
+ f_\theta(\lambda^{-1}\alpha_1v_1)+f_\theta(\lambda^{-1}\beta_1iv_1)+f_\theta(\lambda \alpha_2v_2)&=&0.
\end{eqnarray*}
After dividing by $\lambda^{\delta}$, we deduce that if $\alpha_1,\beta_1,\alpha_2\in\R\backslash\{0\}$, then
\begin{eqnarray}
\label{fthetaj2m-1r1}
 f_\theta(\alpha_1v_1 +\beta_1iv_1+\alpha_2 v_2)- \alpha_1f_\theta(v_1)- \beta_1f_\theta(iv_1)-
 \alpha_2f_\theta(v_2)-&&\\
\label{fthetaj2m-1r2}
- f_\theta(\lambda^{-\delta-1}\alpha_1v_1 +\lambda^{-\delta-1}\beta_1iv_1+\lambda^{1-\delta}\alpha_2 v_2)+&&\\
\label{fthetaj2m-1r3}
+ f_\theta(\lambda^{-\delta-1}\alpha_1v_1)+f_\theta(\lambda^{-\delta-1}\beta_1iv_1)+
f_\theta(\lambda^{1-\delta} \alpha_2v_2)&=&0.
\end{eqnarray}

\noindent{\bf Case 1:} \emph{$f_\theta$ is linear if $\delta\geq 1$.}

Let $\lambda$ tend to infinity. If $\delta=1$, then (\ref{fthetaj2m-1r2}) tends to $-f_\theta(\alpha_2 v_2)$,
and (\ref{fthetaj2m-1r3}) tends to $f_\theta(\alpha_2 v_2)$, yielding (\ref{j2m-1flr}). If $\delta>1$, then both (\ref{fthetaj2m-1r2}) and
(\ref{fthetaj2m-1r3}) tend to zero, and hence 
we conclude (\ref{j2m-1flr}).\\

\noindent{\bf Case 2:} \emph{$f_\theta$ is linear if $-1<\delta<1$ (equivalently, $\delta=0$).}

Since $f_\theta$ is differentiable at $\alpha_2 v_2$ by $\alpha_2>0$ by \eqref{fthetaissmooth}, there exists some $\Omega(\lambda)\in\R$
such that
\begin{eqnarray*}
 f_\theta(\lambda^{-2}\alpha_1v_1 +\lambda^{-2}\beta_1iv_1+\alpha_2 v_2)&=&
 f_\theta(\alpha_2 v_2)+\Omega(\lambda)\\
 \lim_{\lambda\to\infty}\lambda^{1-\delta}\Omega(\lambda)&=&0.
\end{eqnarray*}
Letting $\lambda$ tend to infinity, we conclude  from the $1$-homogeneity of $f_\theta$ that
\begin{eqnarray*}
-f_\theta(\lambda^{-\delta-1}\alpha_1v_1 +\lambda^{-\delta-1}\beta_1iv_1+\lambda^{1-\delta}\alpha_2 v_2)+&&
\\+ f_\theta(\lambda^{-\delta-1}\alpha_1v_1)+f_\theta(\lambda^{-\delta-1}\beta_1iv_1)
 +f_\theta(\lambda^{1-\delta} \alpha_2v_2)&&
=\\
-\lambda^{1-\delta}\Omega(\lambda)+ f_\theta(\lambda^{-1}\alpha_1v_1)&+&f_\theta(\lambda^{-1}\beta_1iv_1) 
\end{eqnarray*}
tends to zero, yieding 
(\ref{j2m-1flr}). \\

\noindent{\bf Case 3:} \emph{$f_\theta$ is linear if $\delta< -1$.}

Let $\lambda$ tend to zero. Since both (\ref{fthetaj2m-1r2}) and
(\ref{fthetaj2m-1r3}) tend to zero, we conclude (\ref{j2m-1flr}).\\

\noindent{\bf Case 4:} \emph{$f_\theta$ is linear if $\delta=-1$.} 

We first observe that if $\delta=-1$ and $\alpha_1,\beta_1,\alpha_2\in\R$, then
\begin{equation}
\label{j2m-1partial-linearity}
f_\theta(\alpha_1v_1 +\beta_1iv_1+\alpha_2 v_2)=f_\theta(\alpha_1v_1+\beta_1 iv_1)+\alpha_2f_\theta(v_2)
\end{equation}
Indeed, letting $\lambda$ tending to zero, (\ref{fthetaj2m-1r1}), (\ref{fthetaj2m-1r2}) and
(\ref{fthetaj2m-1r3}) yield (\ref{j2m-1partial-linearity}).

Hence, it remains to prove that for $\alpha_1,\beta_1\in\R$, we have
$$f_\theta(\alpha_1v_1 +\beta_1iv_1)=\alpha_1f_\theta(v_1)+\beta_1 f_\theta(iv_1).$$

In this case, we set $\mu=\alpha_1$, and consider $\varphi\in {\rm SL}(m,\C)$ by 
$\varphi(v_1)=v_1+\mu v_2$ and $\varphi(v_l)=v_l$
 for $l\geq 2$. 
We observe that 
$\varphi|_L=\tilde{\varphi}=\varphi$  in (\ref{palphabetaodd}).
 Since $\varphi(iv_1)=iv_1+\mu iv_2$ and $\varphi(iv_l)=iv_l$
 for $l\geq 2$, we have
$$
\begin{array}{rclcrcl}
\varphi^{-t}(v_1)&=&v_1& \mbox{ \ and \ }&\varphi^{-t}(iv_1)&=&iv_1\\
\varphi^{-t}(v_2)&=&-\mu v_1+v_2&\mbox{ \ and \ }&\varphi^{-t}(iv_2)&=&-\mu iv_1+iv_2\\
\varphi^{-t}(v_l)&=& v_l
&\mbox{ \ and \ }&\varphi^{-t}(iv_l)&=&i v_l \mbox{ \  for $l>2$.}
\end{array}
$$

For $p,q\in\N$, $0\leq p\leq \theta(2)$ and $0\leq q\leq \theta(m+2)$, let 
$\theta_{pq}\in I$ be such that
\begin{eqnarray*}
\theta_{pq}(1)&=&\theta(1)+p \mbox{ \ and \ }\theta_{pq}(m+1)=\theta(m+1)+q,\\
\theta_{pq}(2)&=&\theta(2)-p \mbox{ \ and \ }\theta_{pq}(m+2)=\theta(m+2)-q,\\
\theta_{pq}(l)&=&\theta(l) \mbox{ \ and \ }\theta_{pq}(m+l)=\theta(m+l) \mbox{ \ for $l=3,\ldots,m$ provided $m\geq 3$.}
\end{eqnarray*}
In particular, $\theta=\theta_{00}$. We observe that if $0\leq p\leq \theta(2)$ and $0\leq q\leq \theta(m+2)$, then
$$
\delta_{pq}=\delta+2p+2q=-1+2p+2q,
$$
and the coefficient of
$\odot_{l=1}^{2m} v_l^{\theta(l)}$ in $\varphi\cdot f$ is
$$
\sum_{p=0}^{\theta(2)}\sum_{q=0}^{\theta(m+2)}
{\theta(1)+p \choose p}{\theta(m+1)+q \choose q}\mu^{p+q}f_{\theta_{pq}}.
$$
We deduce from Case 1 that $f_{\theta_{pq}}$ is $\R$-linear on $\R^{2m}$ unless $p=q=0$.
Therefore it follows from (\ref{palphabetaodd}) applied to the linear combination
$\alpha_1v_1 +\beta_1iv_1+v_2$ that
\begin{eqnarray*}
 f_\theta(\alpha_1v_1 +\beta_1iv_1+v_2)-
\alpha_1 f_\theta(v_1)- \beta_1f_\theta(iv_1)-f_\theta(v_2)-&&\\
- f_\theta((\alpha_1-\mu)v_1 +\beta_1iv_1+ v_2)+
 f_\theta(\alpha_1v_1)+f_\theta(\beta_1iv_1)+f_\theta(-\mu v_1+v_2)&=&0.
\end{eqnarray*}
Since $f_\theta(\alpha_1v_1 +\beta_1iv_1+v_2)=f_\theta(\alpha_1v_1 +\beta_1iv_1)+f_\theta(v_2)$
by \eqref{j2m-1partial-linearity}, and
$$
- f_\theta((\alpha_1-\mu)v_1 +\beta_1iv_1+ v_2)+
 f_\theta(\alpha_1v_1)+f_\theta(\beta_1iv_1)+f_\theta(-\mu v_1+v_2)=0
$$
by $\mu=\alpha_1$ and \eqref{j2m-1partial-linearity}, we conclude Case~4 and the proof of Lemma~\ref{Z2m-1r}.
\proofbox

\subsection{Case $1\leq j\leq m-1$}

The goal of this section is to prove

\begin{lemma}
\label{Zjleqjleqm-1r}
If   $j=1,\ldots,m-1$, then
$Z_j$ is constant zero.
\end{lemma}
\proof 
According to Lemma~\ref{m-1to2m-1odd}, it is sufficient to show for any $L\in{\rm Gr}_{j+1}(\R^{2m})$
of maximal complex rank $j+1$, if $x,y\in L$, and $f:\,L\to \T^r(\R^{2m})$ is the function of Lemma~\ref{m-1to2m-1odd}, then $f(x+y)=f(x)+f(y)$.
There  exists some real orthonormal basis $v_1,\ldots,v_{j+1}$  of  $L$ such that
$x,y\in {\rm lin}_\R\{v_1,v_2\}$, and $v_1,\ldots,v_{j+1}$ can be extended to a complex basis 
$v_1,\ldots,v_m$  of  $\C^m=\R^{2m}$. Therefore it is sufficient to prove that if
 for any $\alpha_1,\alpha_2\in \R$, we have
\begin{equation}
\label{jleqmalphabetar}
f(\alpha_1v_1+\alpha_2v_2)=
\alpha_1f(v_1)+\alpha_2f(v_2).
\end{equation}
As $f$ is continuos, we may assume that
$\alpha_1,\alpha_2\in \R\backslash\{0\}$ in \eqref{jleqmalphabetar}.
We also note that $f$ is $C^\infty$ on $L\backslash \{o\}$ by \eqref{fissmooth}.

To prove \eqref{jleqmalphabetar}, we use the notation leading up to \eqref{palphabetaodd}.
For $\lambda>0$, we define $\varphi\in {\rm SL}(m,\C)$ by $\varphi(v_1)=\lambda v_1$, $\varphi(v_2)=\lambda^{-1} v_2$
and $\varphi(v_l)=v_l$ for $l>2$ where we do not signal the dependence on $\lambda$. It follows that
\begin{eqnarray*}
\varphi^{-t}(v_1)&=&\lambda^{-1}v_1\\
\varphi^{-t}(v_2)&=&\lambda v_2, 
\mbox{ \ moreover  $\varphi^{-t}(v_l)=v_l$ for $l>2$.}
\end{eqnarray*}
For $\theta\in I$, writing
$$
\delta=\theta(1)+\theta(m+1)-\theta(2)-\theta(m+2),
$$
 \eqref{palphabetaodd} yields  
\begin{eqnarray*}
\lambda^{\delta} f_\theta(\alpha_1v_1 +\alpha_2 v_2)-
\lambda^{\delta} \alpha_1f_\theta(v_1)-
\lambda^{\delta} \alpha_2f_\theta(v_2)-&&\\
- f_\theta(\lambda^{-1}\alpha_1v_1 +\lambda\alpha_2 v_2)+
 \alpha_1f_\theta(\lambda^{-1}v_1)+\alpha_2f_\theta(\lambda v_2)&=&0.
\end{eqnarray*}
After dividing by $\lambda^{\delta}$, we deduce that if $\alpha_1,\alpha_2\in\R\backslash 0$, then
\begin{eqnarray}
\label{fthetajleqmr}\nonumber
 f_\theta(\alpha_1v_1 +\alpha_2 v_2)- \alpha_1f_\theta(v_1)- 
 \alpha_2f_\theta(v_2)-&&\\
-f_\theta(\lambda^{-\delta-1}\alpha_1v_1 +\lambda^{1-\delta}\alpha_2 v_2)
+f_\theta(\lambda^{-\delta-1} \alpha_1v_1)+
f_\theta(\lambda^{1-\delta} \alpha_2v_2)&=&0.\qquad
\end{eqnarray}

\noindent{\bf Case 1:} \emph{$f_\theta$ is linear if $\delta\geq 1$.}

If $\lambda$ tends to infinity, then 
$f_\theta(\lambda^{-\delta-1}\alpha_1v_1 +\lambda^{1-\delta}\alpha_2 v_2)$ tends to
$f_\theta(\lambda^{\max\{0,2-\delta\}}\alpha_2 v_2)$,  thus
(\ref{fthetajleqmr}) tends to zero, and in turn
we conclude (\ref{jleqmalphabetar}).\\

\noindent{\bf Case 2:} \emph{$f_\theta$ is linear if $-1<\delta<1$  (equivalently, $\delta=0$).}

Since $f_\theta$ is differentiable at $\alpha_2 v_2$ by $\alpha_2>0$, there exists some $\Omega(\lambda)\in\R$
such that
\begin{eqnarray*}
 f_\theta(\lambda^{-2}\alpha_1v_1 +\alpha_2 v_2)&=&
 f_\theta(\alpha_2 v_2)+\Omega(\lambda),\\
 \lim_{\lambda\to\infty}\lambda^{1-\delta}\Omega(\lambda)&=&0.
\end{eqnarray*}
Letting $\lambda$ tend to infinity in (\ref{fthetajleqmr}), we conclude  from the $1$-homogeneity of $f_\theta$ that
\begin{eqnarray*}
\lambda^{1-\delta}\left(- f_\theta(\lambda^{-2}\alpha_1v_1 +\alpha_2 v_2)+f_\theta( \alpha_2v_2)\right) +\lambda^{-1-\delta}f_\theta(\alpha_1v_1)
&=&\\
-\lambda^{1-\delta}\Omega(\lambda)+ \lambda^{-1-\delta}f_\theta(\alpha_1v_1)&&
\end{eqnarray*}
tends to zero, yielding 
(\ref{jleqmalphabetar}). \\

\noindent{\bf Case 3:} \emph{$f_\theta$ is linear if $\delta\leq -1$.}

If $\lambda$ tends to zero, then (\ref{fthetajleqmr})  tends to zero, and hence 
we conclude (\ref{jleqmalphabetar}).
\proofbox

\medskip
We end this paper by summarizing up the results obtained for even and odd valuations, to give a proof of Proposition~\ref{sl(m,C)invariant}.\\

\noindent{\bf Proof of Proposition~\ref{sl(m,C)invariant}:}\label{finalproof}
Let $Z:\K^{2m}\to\T^r(\R^{2m})$ be an $\SL(m,\C)$-equivariant and translation invariant continuous valuation. McMullen's decomposition Theorem~\ref{McMullendec} yields that $Z$ can be written as
$Z=\sum_{j=0}^{2m}Z_j$ where $Z_j$ is an $\SL(m,\C)$-equivariant and translation invariant continuous $j$-homogeneous valuation for $j=0,\ldots,2m$.
If in addition, $Z$ is even (resp.~odd), then $Z_j$, $j=0,\ldots,2m$, is also even (resp.~odd). 

By Lemmas~\ref{m+1j2m-1}, \ref{Zj1m-1evenr}, \ref{Zmevenr} and  \ref{Z02mevenr}, we have that if $Z$ is even, then each $Z_j$ in McMullen's decomposition is constant zero, and
we conclude Proposition~\ref{sl(m,C)invariant}  for even valuations.

If $Z$ is odd, then combining \eqref{Zoddj02m} and Lemmas~\ref{Zmleqjleq2m-2r}, \ref{Z2m-1r} and \ref{Zjleqjleqm-1r} shows that each $Z_j$ is constant zero. Hence, Proposition~\ref{sl(m,C)invariant} is also proved for odd valuations.
\proofbox

\medskip
\noindent{\bf Acknowledgement: } We are grateful for fruitful discussions with Monika Ludwig, whose ideas opened up new directions while working on this paper, and with Semyon Alesker and Andreas Bernig, who assisted us in some of the known results about translation invariant valuations given in Section 3.

\end{document}